\documentclass{amsart}

\usepackage{amssymb,latexsym}

\numberwithin{equation}{section}

\newtheorem{theorem}{Theorem}[section]

\newtheorem{proposition}[theorem]{Proposition}

\newtheorem{corollary}[theorem]{Corollary}

\newtheorem{definition}[theorem]{Definition}

\newtheorem{lemma}[theorem]{Lemma}

\newtheorem{example}[theorem]{Example}

\newcommand{\credit}[1]{\smallskip\noindent {\textbf{#1.\ }}}

\def\FF{\mathbb{F}}

\def\ii{\mathbf{i}}

\def\RR{\mathbb{R}}

\def\ZZ{\mathbb{Z}}

\begin{document}

\title{Orbits of groups generated by transvections over $\FF_2$}

\author{Ahmet I. Seven}

\email{aseven@lynx.neu.edu}

\thanks{The author's research was supported in part
by Andrei Zelevinsky's NSF grant \#DMS-0200299.}

\date{December 5, 2003}

\maketitle

\section{abstract}

\label{sec:abstract}

Let $V$ be a finite dimensional vector space over the two element field. 
We compute orbits for the linear action of groups generated by
transvections with respect to a certain class of bilinear forms on $V$. 
In particular, we compute orbits that are in bijection with connected components
of real double Bruhat cells in semisimple groups, extending results of
M.Gekhtman, B. Shapiro, M. Shapiro, A.Vainshtein and A. Zelevinsky.

\section{Introduction}

\label{sec:introduction}

Let $V$ be a finite dimensional vector space over the two element field
$\FF_2$ with an $\FF_2$-valued bilinear form $\Omega(u,v)$. 
For any non-zero vector $a\in V$ such that $\Omega(a,a)=0$, the transvection $\tau_a$ is the 
linear transformation  defined as $\tau_a(x)=x+\Omega(x,a)a$ for all $x\in V$.
If the form $\Omega$ is alternating, i.e. $\Omega(a,a)=0$ for all
$a\in V$, then $\tau_a$ preserves $\Omega$,
i.e. $\Omega(x,y)=\Omega(\tau_a(x),\tau_a(y))$; 
in this case $\tau_a$ is called a \emph{symplectic} transvection. 
Since we work over $\FF_2$, each  transvection $\tau_a$ is an involution, i.e.
$\tau_a^2(x)=x$. For a linearly independent subset $B$ of $V$, we denote by $\Gamma_B$ the group 
generated by transvections $\tau_b$ for $b\in B$.  
We define $Gr(B)$ as the graph whose vertex set is $B$ and 
$b_i,b_j$ in $B$ are connected if $\Omega(b_i,b_j)=1$ or $\Omega(b_j,b_i)=1$.
In this paper, we study the linear action of $\Gamma_B$ 
in $V$ for a linearly independent subset $B$ such that $Gr(B)$ is connected. 
For $\Omega$ which is alternating, a description of $\Gamma_B$-orbits are
obtained in \cite{BHII,J,SSVZ} for $B$ such that $Gr(B)$ contains the Dynkin graph $E_6$ as a subgraph
(see Fig.~\ref{fig:E6}).
We give a description of the $\Gamma_B$-orbits for the remaining linearly independent subsets 
(Theorems~\ref{th:D-V000},~\ref{th:xxx},~\ref{non-trivial orbits}). 
Furthermore, we compute $\Gamma_B$-orbits corresponding to a certain class
of non-skew-symmetric bilinear forms (Theorem~\ref{th:GSV}) extending the results of \cite{GSV}. 

\begin{figure}[ht]
\setlength{\unitlength}{1.8pt}

\begin{center}
\begin{picture}(80,20)(0,0)
\thicklines
\put(0,0){\line(1,0){80}}
\put(40,0){\line(0,1){20}}

\put(0,0){\circle*{2.0}$^{x_1}$}
\put(20,0){\circle*{2.0}$^{x_2}$}
\put(40,0){\circle*{2.0}$^{x_3}$}
\put(60,0){\circle*{2.0}$^{x_4}$}
\put(80,0){\circle*{2.0}$^{x_5}$}
\put(40,20){\circle*{2.0}$_{x_6}$}

\end{picture}
\end{center}

\caption{The Dynkin graph $E_6$}
\label{fig:E6}
\end{figure}

Our interest in groups $\Gamma_B$ and their orbits comes from the study
of \emph{double Bruhat cells} initiated in \cite{FZ}. A double Bruhat cell in a 
simply connected connected complex semisimple group $G$ is the variety
$G^{u,v}=BuB\cap B_{-}vB_{-}$, where $B$ and $B_{-}$ are two opposite 
Borel subgroups, and $u$, $v$ any two elements of the Weyl group $W$.
They provide a geometric framework for the study of total positivity
in semisimple groups. They are also closely related to symplectic leaves
in the corresponding Poisson-Lie groups, see e.g. \cite{HKKR}, \cite{KZ}
and references therein. A reduced
double Bruhat cell $L^{u,v}$ is the quotient $G^{u,v}/H$ under the
right or left action of the maximal torus $H=B\cap B_{-}$. It was
shown in \cite{SSVZ} and \cite{Z} that the connected components of the real part
$L^{u,v}(\RR)$ are in a natural bijection with the $\Gamma_{B(\ii)}$-orbits
in $\FF_2^m$, where $\ii$ is a reduced word (of length $m=l(u)+l(v)$) for
the pair $(u,v)$ in the Coxeter group $W\times W$, and $B(\ii)$ is the corresponding set of   
$\ii$-bounded indices as defined in \cite{SSVZ}.
The $\Gamma_{B(\ii)}$-orbits 
for simply laced (resp. non-simply-laced) groups have been computed 
in \cite{SSVZ} (resp. in \cite{GSV}) under the assumption that $Gr(B(\ii))$ contains
the graph $E_6$ (Fig.~\ref{fig:E6}). 
The results presented below are general enough to compute
$\Gamma_{B(\ii)}$-orbits that are related to real double Bruhat cells in semisimple groups.

The groups $\Gamma_B$ also appeared earlier in singularity theory.
To be more precise, let us assume that 
$B$ is a basis and $\Omega$ is an alternating form on $V$.
Our connectedness assumption on $Gr(B)$ implies, in
particular, that $B$ is contained in a $\Gamma_B$-orbit, which
is denoted by $\Delta$. 
In the language of singularity theory,
the orbit $\Delta$ is called a skew-symmetric vanishing lattice with monodromy 
group $\Gamma_B$, c.f. \cite{J}. The main example of a skew-symmetric 
vanishing lattice is the Milnor lattice of an odd dimensional
isolated complete intersection singularity, see e.g. \cite{E}.
A classification of monodromy groups of such lattices is given 
in \cite{J}. According to this classification, 
if the graph $Gr(B)$ contains $E_6$ as a subgraph,
then the group $\Gamma_B$ has precisely two non-trivial $\Gamma_B$-orbits which are the sets $Q_B^{-1}(1)$ and $Q_B^{-1}(0)$, here $Q_B$ is the associated quadratic form \cite{J}. 
To extend this results to an arbitrary basis $B$ (not necessarily containing $E_6$), 
we introduce a function $d:V-\{0\}\to \ZZ_{>0}$ given by
$$d(x)=\min \{s: x=x_1+...+x_s, \mathrm{\:\:for \:\:some\:\:} x_i\in \Delta 
\mathrm{\:\:such\:\: that\:\:} \Omega(x_i,x_j)=0\}.$$
We prove that, for arbitrary $B$, non-trivial $\Gamma_B$-orbits are 
precisely the level sets of $d$ (Theorem~\ref{th:delta decomposition}).
We also give an explicit realization of the function $d$ in terms of the graph $Gr(B)$
for any basis $B$ (Theorems~\ref{th:delta decomposition}, \ref{th:xxx}). 
Furthermore we extend this realization to a linearly independent set $B$ which is not a basis
and give an explicit description of the orbits (Theorem~\ref{non-trivial orbits}) 
extending the results of \cite{SSVZ}.





To study the action of $\Gamma_B$, we use combinatorial and algebraic
methods. Our main combinatorial tool is a class of graph transformations 
generated by \emph{basic moves}. More precisely, for every two elements $a,c \in B$
such that $\Omega(a,c)=1$, the basic move $\phi_{c,a}$ 
replaces $c$ with $\tau_a(c)$ and leaves other elements of $B$ unchanged.
The essential feature of those moves is to preserve the associated group $\Gamma_B$,
i.e. if $B'$ is obtained from $B$ by a sequence of basic moves, then
$\Gamma_B=\Gamma_{B'}$. Basic moves were suggested to me by A. Zelevinsky; however 
it was brought to my attention that they had been introduced in \cite{BHI}.
(We thank the anonymous referee for pointing this out).
It is important, e.g. in the theory of double Bruhat cells, to be able to
recognize whether a given graph can be obtained from another using basic moves. 
Our Theorems~\ref{th:D-V000}, \ref{th:D-V000-char}, \ref{th:D-V000-span}, \ref{th:xxx} 
solve this \emph{recognition} problem for the classes of graphs that do not contain the subgraph $E_6$.

\section{Main Results}

\label{sec:background}

In this section, we recall some statements
from \cite{BHI,BHII,J} and state our main results.
As in Section~\ref{sec:introduction}, $V$ is
a finite dimensional vector space over the 2-element field $\FF_2$.
Unless otherwise stated, $\Omega$ denotes an $\FF_2$-valued alternating bilinear form on $V$. For a subspace $U$ of $V$, we denote by $U_0$ the
kernel of the form $\Omega$ in $U$, i.e.
$U_0=\{x \in U: \Omega(x,u)=0, \: \mathrm{for \: any} \: u\in  U\}$.
We assume that $B$ is a linearly independent subset whose graph $Gr(B)$ is connected. By some abuse of notation, we will sometimes denote $Gr(B)$ by $B$. By a \emph{subgraph} of $Gr(B)$, we always mean a graph $X$ 
obtained from $Gr(B)$ by taking a subgraph on a subset of
vertices. 
If $B$ is a basis, then there is a one-to-one correspondence $x\to Gr(B,x)$ between $V$ and
subgraphs of $Gr(B)$ defined as follows: $Gr(B,x)$ is the subgraph of $Gr(B)$ on vertices 
$b_{i_1},..,b_{i_k}$, where $x=b_{i_1}+...+b_{i_k}$ is the expansion of $x$ in the basis $B$. 
By some abuse of notation, we sometimes denote $Gr(B,x)$ by $x$.
We say that a vector $u$ is contained in $Gr(B,x)$ if
$Gr(B,u)$ is contained in $Gr(B,x)$.

\begin{definition}

\label{def:moves} 
\cite[Definition~3]{BHI}
\rm {Let $a,c$ be in $B$ such that $\Omega(a,c)=1$.
The \emph{basic move} $\phi_{c,a}$ is the transformation that replaces
$c\in B$ by $\tau_a(c)=c+a$ and keeps other elements the same, 
i.e. $\phi_{c,a}(c)=c+a$ and $\phi_{c,a}(b)=b$ for $b\ne c$. We call 
two linearly independent subsets $B$ and $B'$ equivalent
if $B'$ is obtained from $B$ by a sequence of basic moves.
If $B$ and $B'$ are equivalent, then their graphs $Gr(B)$ and $Gr(B')$
are also said to be equivalent.}

\end{definition}

\noindent
We note that this equivalence relation is well-defined because 
$\phi_{\tau_a(c),a}\phi_{c,a}(B)=B$. We also note that  
the basic move $\phi_{c,a}$ changes $Gr(B)$ as follows: 
suppose that $p$ and $q$ are the vertices that represent the basis vectors $a$ and $c$ respectively. 
The move $\phi_{c,a}$ 
connects $q$ to vertices that are connected to $p$ but not connected to $q$. 
At the same time, it disconnects vertices from $q$ if they are connected to $p$. 
It follows that, for any basis $B'$ which is equivalent to $B$, 
the graph $Gr(B')$ is connected. Furhermore, $\Gamma_B=\Gamma_{B'}$ \cite[Proposition~3.1]{BHI}.
Also, there exists a basis $B'$ equivalent to 
$B$ such that $Gr(B')$ is a tree \cite[Theorem~3.3]{BHII}.

\begin{figure}[ht]
\setlength{\unitlength}{1.8pt}
\begin{center}
\begin{picture}(70,20)(-10,-10)
\thicklines
\put(-10,0){\circle*{2.0}$^{a_1}$}
\put(0,0){\circle*{2.0}$^{a_2}$}
\put(10,0){\circle*{2.0}}
\put(20,0){\circle*{2.0}}
\put(30,0){\circle*{2.0}}
\put(40,0){\circle*{2.0}}
\put(50,0){\circle*{2.0}$^{\:a_m}$}
\put(60,0){\circle*{1.0}$_{c_i}$}
\put(60,10){\circle*{2.0}$_{c_1}$}
\put(60,-10){\circle*{2.0}$_{c_k}$}
\put(60,5){\circle*{1.0}}
\put(60,-5){\circle*{1.0}}
\put(-10,0){\line(1,0){70}}
\put(50,0){\line(1,1){10}}
\put(50,0){\line(1,-1){10}}

\end{picture}
\end{center}
\caption{The graph $D_{m,k}$}
\label{fig:broom}
\end{figure}
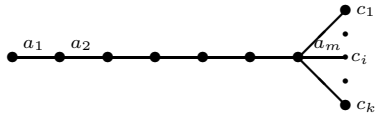

\begin{definition}
\label{def:type} \rm{ Let $m,k$ be integers such that $m \geq 2$, $k \geq 1$.
A graph of the form in Fig.~\ref{fig:broom} is said to be of type $D_{m,k}$.}
\end{definition}

\subsection{Orbits of groups generated by symplectic transvections of a basis}
\label{sub:intro-basis}
To describe such orbits we need to recall some basic facts from the theory of 
quadratic forms over $\FF_2$.
A quadratic form $Q$ is an $\FF_2$-valued function on $V$
having the following property:
$$Q(u+v)=Q(u)+Q(v)+g(u,v), (\mathrm{for\:\: all} \: u,v \in V)$$ 
where $g:V\times V \to \FF_2$ is an alternating bilinear form.
It is clear that the quadratic form $Q$ completely determines the 
associated bilinear form $g$. Recall (see e.g.\cite{D}) that 
there exists a $symplectic$ basis $\{e_1,f_1,...,e_r,f_r,h_1,...,h_p\}$ in $V$ such that
$g(e_i,f_j)=\delta_{i,j}$, 
and the rest of the values of $g$ are $0$;
here $\delta_{i,j}$ is the Kronecker symbol. Let us write 
$V_0=\{x \in V: g(x,v)=0, \: \mathrm{for \: any} \: v\in  V\}$.
If $Q(V_0)=\{0\}$, then the Arf invariant of $Q$ is defined as 
$$Arf(Q)=\sum Q(e_i)Q(f_i)$$
It is well known from the theory of quadratic forms that $Arf(Q)$ is independent of the
choice of the symplectic basis (\cite{D}).

Two quadratic forms $Q$ and $Q'$ on $V$ are 
isomorphic if there is a
linear isomorphism $T:V\to V$ such that $Q(T(x))=Q'(x)$ for any $x\in V$.
According to \cite{C,D}, isomorphism classes of quadratic forms $\{Q\}$ on $V$  
are determined by their Arf invariants and 
their restrictions $\{Q \vert _ {V_0}\}$.
More precisely, for fixed dimensions of $V$ and $V_0$, there exist at most 
3 isomorphism classes of quadratic forms $\{Q\}$ and each isomorphism class is determined 
by precisely one of the following:

\begin{itemize} 
\item[{\rm(i)}] 
 $Q(V_0)=0, Arf(Q)=1$
\item[{\rm(ii)}] 
$Q(V_0)=0, Arf(Q)=0$
\item[{\rm(iii)}] 
$Q(V_0)=\FF_2$
\end{itemize}

Let us now assume that $B$ is a basis of the $\FF_2$-space $V$ equipped with an alternating form $\Omega$. We denote by $Q_B$ the unique quadratic form associated with 
$\Omega$ and $B$ as follows: $Q_B(u+v)=Q_B(u)+Q_B(v)+\Omega(u,v), (u,v \in V)$,
and $Q_B(b)=1$ for all $b \in B$. It is easy to see that $Q_B$ is $\Gamma_B$-invariant, i.e. $Q_B(\tau_a(x))=Q_B(x)$ for all $a\in B$. This also implies that quadratic forms are invariant under basic moves; 
i.e. if $B$ and $B'$ are equivalent bases, then 
$Q_B(x)=Q_{B'}(x)$ for any $x \in V$. Furthermore, the function $Q_B$
completely determines the $\Gamma_B$-orbits in $V-V_0$ when the graph $Gr(B)$ 
contains a subgraph equivalent to $E_6$:

{
\begin{theorem}

\label{th:Janssen classification}

\cite[Theorem~4.1]{BHII}, \cite[Theorem~3.8]{J}
{Let $V$ and $V_0$ have dimensions $2n+p$ and $p$ respectively.
Suppose that $B$ 
is equivalent to a tree which 
contains $E_6$ as a subgraph. Then $B$ is equivalent
to one of the trees in Fig.~\ref{fig:JanssenA}, Fig.~\ref{fig:JanssenB} or
Fig.~\ref{fig:JanssenC}, depending on $Q_B(V_0)$ and $Arf(Q_B)$.
Furthermore, the group $\Gamma_B$ has precisely two orbits in $V-V_0$. 
They are intersections of $V-V_0$ with the sets $Q_B^{-1}(0)$ and $Q_B^{-1}(1)$.

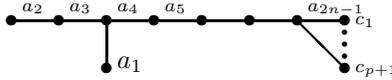
\begin{figure}[ht]
\setlength{\unitlength}{1.8pt}
\begin{center}
\begin{picture}(70,40)(-10,-20)
\thicklines
\put(-10,0){\circle*{2.0} $^{a_2}$}
\put(0,0){\circle*{2.0} $^{a_3}$}
\put(10,0){\circle*{2.0} $^{a_4}$}
\put(10,-10){\circle*{2.0} $a_1$}
\put(20,0){\circle*{2.0} $^{a_5}$}
\put(30,0){\circle*{2.0}}
\put(40,0){\circle*{2.0}}
\put(50,0){\circle*{2.0}$^{a_{2n-1}}$}
\put(60,0){\circle*{2.0} $_{c_1}$}
\put(60,-2.5){\circle*{1.0}}
\put(60,-5){\circle*{1.0}}
\put(60,-7.5){\circle*{1.0}}
\put(60,-10){\circle*{2.0} $_{c_{p+1}}$}
\put(-10,0){\line(1,0){70}}
\put(10,0){\line(0,-1){10}}
\put(50,0){\line(1,-1){10}}

\end{picture}

\end{center}

\caption{$Arf(Q_B)=1$ if $n\equiv 2,3$ mod(4), $Arf(Q_B)=0$ if $n\equiv 0,1$ mod(4)
, $V_0=$ linear span of $\{c_1+c_{p+1},c_2+c_{p+1},...,c_p+c_{p+1}\}$, $Q_B(V_0)$=0.}

\label{fig:JanssenA}

\end{figure}

\begin{figure}[ht]
\setlength{\unitlength}{1.8pt}

\begin{center}
\begin{picture}(70,40)(-10,-20)

\thicklines

\put(-10,0){\circle*{2.0} $^{a_3}$}
\put(0,0){\circle*{2.0} $^{a_4}$}
\put(10,0){\circle*{2.0} $^{a_5}$}
\put(20,-10){\circle*{2.0} $_{a_2}$}
\put(20,-20){\circle*{2.0}$_{a_1}$}
\put(20,0){\circle*{2.0} $^{a_6}$}
\put(30,0){\circle*{2.0}}
\put(40,0){\circle*{2.0}}
\put(50,0){\circle*{2.0}$^{a_{2n-1}}$}
\put(60,0){\circle*{2.0} $_{c_1}$}
\put(60,-2.5){\circle*{1.0}}
\put(60,-5){\circle*{1.0}}
\put(60,-7.5){\circle*{1.0}}
\put(60,-10){\circle*{2.0} $_{c_{p+1}}$}

\put(-10,0){\line(1,0){70}}

\put(20,0){\line(0,-1){20}}

\put(50,0){\line(1,-1){10}}

\end{picture}

\end{center}

\caption{$Arf(Q_B)=1$ if $n\equiv 0,1$ mod(4), $Arf(Q_B)=0$ if $n\equiv 2,3$ mod(4)
, $V_0=$ linear span of $\{c_1+c_{p+1},c_2+c_{p+1},...,c_p+c_{p+1}\}$, $Q_B(V_0)=0$.}
\label{fig:JanssenB}
\end{figure}
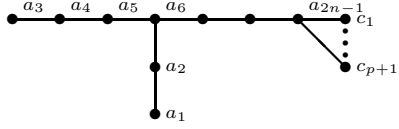

\begin{figure}[ht]
\setlength{\unitlength}{1.8pt}

\begin{center}
\begin{picture}(70,40)(-10,-10)
\thicklines
\put(-10,0){\circle*{2.0}$^{a_2}$}
\put(0,0){\circle*{2.0}$^{a_3}$}
\put(10,0){\circle*{2.0}$^{a_4}$}
\put(20,-10){\circle*{2.0}$_{a_1}$}
\put(20,0){\circle*{2.0}$^{a_5}$}
\put(30,0){\circle*{2.0}}
\put(40,0){\circle*{2.0}}

\put(50,0){\circle*{2.0}$^{a_{2n}}$}
\put(60,0){\circle*{2.0}$_{c_1}$}
\put(60,-2.5){\circle*{1.0}}
\put(60,-5){\circle*{1.0}}
\put(60,-7.5){\circle*{1.0}}
\put(60,-10){\circle*{2.0}$_{c_p}$}

\put(-10,0){\line(1,0){70}}
\put(20,0){\line(0,-1){10}}
\put(50,0){\line(1,-1){10}}

\end{picture}
\end{center}
\caption{$Q_B(V_0)=\FF_2$,
$V_0=$linear span of $\{c_1+c_p,c_2+c_p,...,c_{p-1}+c_p,a_1+a_2+a_4\}$}
\label{fig:JanssenC}
\end{figure}
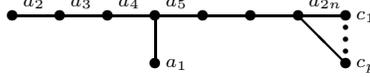
}
\end{theorem}
}

\begin{proposition}
\label{pr:E_6 ne D}
\cite[Theorem~4.2]{BHII}
 Suppose that {\rm{dim($V$)}}$\geq 3$ and let $B$ be a tree that does not contain $E_6$. 
Then it is equivalent to a tree of type $D_{m,k}$, i.e. a tree of the form in Fig.~\ref{fig:broom},
 for some $m\geq 2, k\geq 1$.
Furthermore, $B$ is not equivalent to any tree $B'$ that contains $E_6$.
\end{proposition}

Thus the trees given in Figures~\ref{fig:broom}, \ref{fig:JanssenA},
\ref{fig:JanssenB}, \ref{fig:JanssenC}
form a complete set of representatives for equivalence classes of graphs.
For a basis $B$ which is a tree of type $D_{m,k}$, there is the following
description of $\Gamma_B$-orbits.

\begin{theorem}
\label{pr:broom}
\cite[Theorem~7.2]{BHII}
Let $V$ be an $\FF_2$-space with $dim(V)\geq 3$ and let $B$ be a basis
whose graph is a tree of type $D_{m,k}$ , i.e. a tree of the form in Fig.~\ref{fig:broom},
with $m\geq 2, k\geq 1$.
Then, any two vectors $x, y$ in $V-V_0$ lie in the same $\Gamma_B$ orbit
if and only if $Gr(B,p(x))$ and $Gr(B,p(y))$ have the same number of connected components,
where $p$ is the linear map on $V$ defined as $p(a_i)=a_i$ for all $i=1,...,m$, $p(c_j)=c_1$ for $j=1,...,k$.
\end{theorem}

Our next result extends Theorem~\ref{th:Janssen classification} (and Theorem~\ref{pr:broom}) to an arbitrary basis.


\begin{theorem}
\label{th:delta decomposition}
Let $B$ be an arbitrary basis in $V$ such that $Gr(B)$ is connected
and let $d:V\to \ZZ_{>0}$ be the function defined as 
$$d(x)=\min \{s: x=x_1+...+x_s, \mathrm{\:\:for \:\:some\:\:} x_i\in \Delta 
\mathrm{\:\:such\:\: that\:\:} \Omega(x_i,x_j)=0\}$$
where $\Delta$ is the $\Gamma_B$-orbit that contains $B$.
Let $x,y$ be vectors in $V-V_0$. Then $x$ and $y$ lie 
in the same $\Gamma_B$-orbit if and only if $d(x)=d(y)$.
Furthermore, if $B$ is equivalent to a tree which contains $E_6$, then $d(x)=2-Q_B(x)$ for any $x\in V-V_0$.
\end{theorem}

In Theorem~\ref{th:xxx}, we will obtain an explicit expression of the function $d$ for an arbitrary basis $B$ which
is equivalent to a tree that does not contain $E_6$. Our next result allows one to recognize such bases easily.

\begin{theorem}
\label{th:subE6}
A connected graph $B$ is equivalent to a tree of type $D_{m,k}$
with $m\geq 2, k\geq 1$
if and only if it does not contain any subgraph 
which is equivalent to $E_6$.
Furthermore, the graph $B$ is equivalent to a tree of type $D_{m,1}$
if and only if it does not have subgraphs of the following form:
\begin{align} 
\label{eq:A} 
&\text{a tree of type $D_{2,2}$,}
\\
\nonumber 
&\text{two triangles sharing a common edge,} 
\\
\nonumber 
&\text{any cycle whose length is greater than or equal to 4.} 
\end{align} 
\end{theorem}
\noindent
Corollary~\ref{cor:sub-connected} also provides an alternate method to recognize a graph which is equivalent
to a tree that does not contain $E_6$.

Our next theorem introduces an important class of vectors which are fixed by $\Gamma_B$.

\begin{theorem}
\label{th:V000}

Let $B$ be an arbitrary basis in $V$ such that $Gr(B)$ is connected
and let $\Delta$ denote the $\Gamma_B$-orbit that contains $B$. 
Let $V_{000}$ be the set that consists of vectors $y$ in $V_0$ 
such that $y=x_1+x_2$ for some $x_1, x_2 \in \Delta$.
Then $V_{000}$ is a vector subspace of $V$ and   
every $\Gamma_B$-orbit in $V-V_0$ is a union of cosets in $V/V_{000}$.
Furthermore, we have {\rm{dim($V_0/V_{000}$)}}$\leq 1$.

\end{theorem}

Our next results allow one to locate all of the vectors in $V_{000}$ for a basis
which does not contain $E_6$.

\begin{theorem}
\label{th:D-V000}

Let $B$ be a basis whose graph is equivalent to a tree of type $D_{m,k}$ with $m\geq 2, k\geq 1$. 
Suppose that $X$ is a subgraph which is one of the types in \eqref{eq:A}. Then $\FF_2^X \cap V_{000} \ne \{0\}$, where $\FF_2^X$ is the linear span of the vectors contained in $X$. Furthermore, if $X=[x_1,x_2,...,x_r]$ is a cycle whose length $r$ is greater than or equal to 5, then $\FF_2^X \cap V_{000}$ is spanned by the vector $x_1+x_2+...+x_r$.






\end{theorem}

\begin{theorem}
\label{th:D-V000-char}

The conditions in Theorem~\ref{th:D-V000} characterize a graph which is equivalent to a tree of type $D_{m,k}$ for some $m\geq 3, k\geq 1$.

\end{theorem}

\begin{theorem}
\label{th:D-V000-span}

Let $B$ be a basis whose graph is equivalent to a tree of type $D_{m,k}$ with $m\geq 2, k\geq 1$. If {\rm{dim($V$)}}$\geq 3$, then $V_{000}=\sum_X \FF_2 ^X\cap V_{000}$ where $X$ runs through the subgraphs in \eqref{eq:A}.

\end{theorem}

\begin{corollary}
\label{cor:sub-connected}
Let $B$ be a basis whose graph is 
equivalent to a tree of type $D_{m,k}$ with $m\geq 2, k\geq 2$. 
Suppose that $b$ is an arbitrary vertex contained in $Gr(B,u)$ for some 
$u \in V_{000}$. Then the graph $Gr(B-\{b\})$ is connected and it
is equivalent to a tree of type $D_{m,k-1}$. 
\end{corollary}

Our next result gives an explicit expression for the function $d$ defined in 
Theorem~\ref{th:delta decomposition}. 

\begin{theorem}
\label{th:xxx}

Let $B$ be a basis which is equivalent to
a tree of type $D_{m,k}$ with $m\geq 2, k\geq 1$. 
Suppose that $x\in V-V_0$ and let $\bar{x}$ be any minimal representative in the coset $x+V_{000}$
, i.e. $Gr(B,\bar{x})$ does not contain any non-zero vector in $V_{000}$. Then 
\begin{equation}
\label{eq:xxx}
d(x)=c(\bar{x})+\sum_{A:\vert A \vert \geq 3}(\lceil\: \vert A \vert /2\:\rceil -1)
\end{equation}
where $c(\bar{x})$ is the number of connected components
of $\bar{x}$ and $A$ runs through the set of maximal
complete subgraphs of $\bar{x}$ and $\vert A \vert$ is the
number of vertices in $A$. 
\end{theorem}

\subsection{Orbits of groups generated by symplectic transvections of a linearly independent subset which is not a basis}
\label{sub:intro-non-basis}
We first recall the following statement from \cite{SSVZ}.

\begin{theorem}
\label{zelevinsky}
\cite[Theorem~6.2]{SSVZ}
Let $B$ be a linearly independent subset which is not a basis in $V$
and let $U$ denote the linear span of $B$. Suppose that $Gr(B)$ is connected.
If $Gr(B)$ contains a subgraph which is equivalent to $E_6$, then $\Gamma_B$ has precisely two orbits in $(v+U)-V^{\Gamma_B}$, where $V^{\Gamma_B}$ is the set of vectors
which are fixed by $\Gamma_B$. They are intersections of $(v+U)-V^{\Gamma_B}$ with the sets 
$Q_{B\cup\{v\}}^{-1}(0)$ and $Q_{B\cup\{v\}}^{-1}(1)$.
\end{theorem}

Our next result is a counterpart of Theorem~\ref{zelevinsky} for the graphs not containing $E_6$.

\begin{theorem}
\label{non-trivial orbits}
Let $B$ be a linearly independent subset which is not a basis in $V$.
Suppose that $Gr(B)$ is connected and does not contain a subgraph which is equivalent to $E_6$. Suppose also that $dim(U) \geq 2$ where $U$ denotes the linear span of $B$. 
Then we have the following statements:

\begin{itemize}
\item[{\rm(i)}] 
If $\Omega(v,U_{000})\ne \{0\}$, then there are precisely two
$\Gamma_B$-orbits in $v+U$. They are
are intersections of $v+U$ with the sets $Q_{B\cup\{v\}}^{-1}(0)$ and 
$Q_{B\cup\{v\}}^{-1}(1)$.

\item[{\rm(ii)}] 
If $\Omega(v,U_{000})=\{0\}$ 
then $\Gamma_B$-orbits in $v+U$ are intersections of
$v+U$ with the $\Gamma_{B\cup \{w\}}$-orbits in the linear span of $B\cup \{w\}$ for some $w \in v+U$ such that $Gr(B\cup\{w\})$ does not contain any subgraph which is equivalent to $E_6$.

\end{itemize}
\end{theorem}

\subsection{Orbits of groups generated by non-symplectic transvections}
\label{sub:intro-non-skew}
Let us first recall that, for an arbitrary bilinear form $\Omega$ on $V$
and a linearly independent set $B$, 
we define its graph $Gr(B)$ as the graph whose vertex set is $B$ and 
$b_i,b_j$ in $B$ are connected if $\Omega(b_i,b_j)=1$ or $\Omega(b_j,b_i)=1$.
The following result is a generalization of \cite[Theorem~3]{GSV}:

\begin{theorem}
\label{th:GSV}

Suppose that $\Omega$ is a non-alternating bilinear form on the $\FF_2$-space $V$.
Let $B$ be a linearly independent subset whose graph
is connected. 
Suppose that there exists a disjoint collection of connected graphs $B_1,...,B_r$ with $B=\cup B_i$ such that 
$\Omega \vert_{\FF_2^{B_i}}$ is alternating and, for any $b_i \in B_i, b_j\in B_j$ with  $i\ne j$, we have $\Omega(b_i,b_j)=1$ if and only if $j=i+1$. Suppose that $x\in V$ is not fixed by $\Gamma_B$ and let 
$L=min\{j:x \mathrm{\:is\: not\: fixed\: by\:} \Gamma_{B_j}\}$. Then the $\Gamma_B$-orbit of $x$ coincides with the set $\Gamma_{B_L}(\pi_L(x))+\FF_2^{\cup_{j>L}B_j}$ where $\Gamma_{B_L}(\pi_L(x))$ is the $\Gamma_{B_L}$-orbit
of the vector $\pi_L(x)$; here $\pi_L$ is the standard projection onto $\FF_2^{B_L}$.  
\end{theorem}


\section{Orbits of groups generated by symplectic transvections of a basis}

\label{sec:basis}

In this section, we prove 
Theorems~\ref{th:delta decomposition}--\ref{th:xxx}. 
It will be convinent for us to prove, first, Theorem~\ref{th:V000}.
Let us keep the notation of Section~\ref{sec:background}
and recall that $V_{000}$ is the set  that consists of vectors $y$ in $V_0$ such that $y=x_1+x_2$ for some $x_1, x_2 \in \Delta$. 
We also introduce the subspace $V_{00}=\{y \in V_0: Q_B(y)=0\}$
where $Q_B$ is the (unique) quadratic function associated with $B$
(c.f. Section~\ref{sec:background}). We note that 
$V_{000}\subset V_{00}$.

\subsection{Proof of Theorem~\ref{th:V000}}

To prove the theorem, it is enough to show that for any $x \in V-V_0$ and $u \in V_{000}$, 
the vectors $x$ and $x+u$ lie in the same $\Gamma_B$-orbit. 
Let us assume that $u=u_1+u_2$, where $u_1, u_2\in \Delta$. Since $B$ lies in $\Delta$, we may also assume, without loss of generality that, $u_1 \in B$. We note that $\Omega(u_1,u_2)=0$
since $u_1+u_2 \in V_0$.

We claim that there exists $\gamma\in \Gamma_B$ such that $\Omega(\gamma(x),u_1)=1$.
Suppose that $\Omega(\gamma(x),u_1)=0$ for all $\gamma\in \Gamma_B$. 
Then $\Omega(x,\gamma (u_1))=0$ for all $\gamma\in \Gamma_B$.
This implies, in particular, that $\Omega(x,b)=0$ for all $b\in B$, because $u_1\in B$ lies in the orbit $\Delta$. 
This would imply that $x\in V_0$, resulting in a contradiction because $x\notin V_0$.
We note that $\Omega(\gamma(x),u_2)=1$ because $u=u_1+u_2\in V_0$

Let us now consider
$\gamma^{-1}\tau_{u_2}\tau_{u_1}\gamma$. This automorphism is in $\Gamma_B$ because $\tau_{u_2}\in \Gamma_B$ by 
\cite[Proposition~3.1]{BHI}. Then $\gamma^{-1}\tau_{u_2}\tau_{u_1}\gamma(x)=x+u_1+u_2=x+u$, hence $x$ and $x+u$ are contained in the same $\Gamma_B$-orbit. 

The remaining part of the Theorem~\ref{th:V000} follows from the following statement.

\begin{proposition}
\label{th:complement of U000}

The spaces $V_0,V_{00}$ and $V_{000}$
have the following descriptions in terms of the trees in 
Figures~\ref{fig:broom},\ref{fig:JanssenA},\ref{fig:JanssenB},\ref{fig:JanssenC}.
\begin{itemize} 
\item[{\rm(i)}] 
Suppose that $B$ is of type $D_{m,k}$ with $m\geq 2, k \geq 1$
indexed as in Fig.~\ref{fig:broom}.

If $m$ is odd, then $V_0=V_{00}=V_{000}=${\rm linear span of} 
$\{c_1+c_2, c_1+c_3,...,c_1+c_k\}$. 

If $m=2$, then $V_0=V_{00}=V_{000}=${\rm linear span of} 
$\{a_1+c_1,c_1+c_2, c_1+c_3,...,c_1+c_k\}$.

If $m>2$ and $m\equiv 2$ $mod(4)$, then $V_{000}=${\rm linear span of }
$\{c_1+c_2, c_1+c_3,...,c_1+c_k\}$ and $V_0=V_{00}=${\rm linear span of} 
$V_{000}\cup \{a_1+a_3+...+a_{m-1}+c_1\}$.  

If $m>2$ and $m\equiv 0$ $mod(4)$, then $V_{00}=V_{000}=${\rm linear span of }
$\{c_1+c_2, c_1+c_3,...,c_1+c_k\}$ and $V_0=${\rm linear span of }
($V_{000}\cup \{a_1+a_3+...+a_{m-1}+c_1\}$). 

\item[{\rm(ii)}] 
If $B$ is as in Fig.~\ref{fig:JanssenA} or Fig.~\ref{fig:JanssenB}, then
$V_0=V_{00}=V_{000}=${\rm linear span of} $\{c_1+c_2, c_1+c_3,...,c_1+c_k\}$ 

\item[{\rm(iii)}] 
If $B$ is as in Fig.~\ref{fig:JanssenC}, 
then $V_{00}=V_{000}=${\rm linear span of} $\{c_1+c_2, c_1+c_3,...,c_1+c_k\}$ and
$V_0=${\rm linear span of} $V_{000}\cup \{a_1+a_2+a_4\}$.
\end{itemize} 

\end{proposition}








\begin{corollary}
\label{cor:mk}

For any basis $B$ whose graph is equivalent to a tree of type $D_{m,k}$
with $m\geq 2, k\geq 1$, we have the following:

\begin{itemize} 
\item[{\rm(i)}] 

If $m > 2$, then $dim(V_{000})=k-1$. If $m=2$, then $dim(V_0)=dim(V_{000})=k$.

\item[{\rm(ii)}] 

The number of fixed points of $\Gamma_B$ in $V$ is $2^{k-1}$ if 
$m$ is odd, and $2^k$ if $m$ is even.

\item[{\rm(iii)}] 
The number of $\Gamma_B$-orbits which are not fixed points 
is $(m+1)/2$ if 
$m$ is odd, and $m/2$ if $m$ is even.



\end{itemize} 
\end{corollary}

\subsection{Proof of Theorem~\ref{th:delta decomposition}}
\label{sec:delta decomposition}

Let us first introduce some terminology.
For a vector $x \in V$, we say that 
$x=x_1+...+x_s$ is a $\Delta$-decomposition of $x$ if $x_1,...,x_s\in \Delta$
and $\Omega(x_i,x_j)=0$ for all $i,j=1,...,s$ 
, here $\Delta$ denotes the $\Gamma_B$-orbit that contains $B$ (and any basis equivalent to it). If $Q_B(x)=1$ (resp.
$Q_B(x)=0$), then any $\Delta$-decomposition of $x$ has an odd (resp. even) number of components. Let us also note that 
$$d(x)=\min \{s: x \: \mathrm{has \: a} \: \Delta \: \mathrm{decomposition}\: 
x=x_1+...+x_s\}.$$

We first prove the existence of a $\Delta$-decomposition for an arbitrary $x\in V-V_0$. Since any connected graph is equivalent to a tree, we may assume that $B$ is a tree. 
Let $B_1,..., B_s$ be the connected components of $Gr(B,x)$ and let
$x_1,...,x_s$ be the corresponding vectors in $V$. We claim that
$x=x_1+....+x_s$ is a $\Delta$-decomposition for $x$. We note
that $\Omega(x_i,x_j)=0$ for all $i\ne j$. To show that $x_i\in \Delta$
for $i=1,...,s$, it is enough to show, without loss of generality,
that $x_1 \in \Delta$. Let us assume $x_1=b_1+...+b_k, b_i\in B$.
Since $Gr(B,x_1)$ is a tree, it has a leaf, i.e. a vertex, say $b_k$,
which is connected to precisely one vertex in $Gr(B,x)$.Then 
$Gr(B,\tau_{b_k}(x))$ is a tree with $k-1$ vertices. By induction,
we obtain $\gamma \in \Gamma_B$ such that $\gamma(x_1) \in B$, i.e. $x_1\in \Delta$.  

Let us also note that, for any $\Delta$-decomposition $x=x_1+...+x_s$ and any $\gamma \in \Gamma$, we have a $\Delta$-composition $\gamma(x)=\gamma(x_1)+...+\gamma(x_s)$. Thus, we have 

\begin{align}
\label{eq:d-inv}
&\text{$d(x)=d(\gamma(x))$ for any $\gamma \in \Gamma.$}
\end{align}

\begin{lemma}
\label{}
Suppose that $B$ is equivalent to a tree that contains $E_6$. Then, for any $x\in V-V_0$,

\begin{itemize} 
\item[{\rm(i)}] 
$d(x)=1$ if and only if $Q_B(x)=1$
\item[{\rm(ii)}] 
$d(x)=2$ if and only if $Q_B(x)=0$
\end{itemize}

\end{lemma}

The "only if" parts are clear. To prove the "if" parts, let us first assume $Q_B(x)=1$.
Then $x\in \Delta$, by Theorem~\ref{th:Janssen classification}, hence $d(x)=1$. Let us now assume that $Q_B(x)=0$ and let $x=x_1+...+x_{2l}$, be a $\Delta$-decomposition. 
Then there exist $j \in \{1,...,2l\}$ such that $x+x_j \notin V_0$ (otherwise
$x=\sum(x+x_j) \in V_0)$. Since $x=(x+x_j)+x_j$,
and $Q_B(x+x_j)=Q_B(x_{j})=1$, we have $d(x)=2$ and we are done.

To complete the proof of the theorem let us now assume that $B$ is equivalent to a tree $B'$ which is of type $D_{m,k}$.
We claim that, for any $x\in V-V_0$, the number $d(x)$ is equal to the number of connected components of $Gr(B',p(x))$
, where $p$ is the function defined in Theorem~\ref{pr:broom}. In view of Theorem~\ref{th:V000} and \eqref{eq:d-inv}, this claim follows from the following stronger statement.

\begin{lemma}
\label{lem:A-d}

Let $B$ be a basis equivalent to a tree $B'$ of type $D_{m,k}$ with $m\geq 2$ and $k\geq 1$. Suppose that $x\in V$ and let $x=x_1+...+x_d$ be a $\Delta$-decomposition with $d=d(x)$ such that, 
\begin{align}
\label{eq:A-d}
&\text{for any $i\ne j$, we have $x_i+x_j\notin V_{000}$.} 
\end{align}
Then there exists $\alpha \in \Gamma_B$ such that $\alpha(x)=y_1+y_3+...+y_d$, where $y_i=\alpha(x_i) \in B'$.
\end{lemma}

Let us assume that $B'$ is indexed as in Fig.~\ref{fig:broom}.
Applying if necessary an element of $\Gamma_B$, we may assume that $x_1=a_1$.
By \eqref{eq:A-d}, all $x_i$ for $i\geq 2$ are 
contained in $\{a_3,...,a_m,c_i\}$. Our claim follows by induction on $d$.
This completes the proof of Theorem~\ref{th:delta decomposition}.

\subsection{Proof of Theorem~\ref{th:subE6}}
\label{subsec:th-E6}


We will first prove that
\begin{align} 
\label{eq:A2} 
&\text{a graph $B$ is equivalent to a tree of type $D_{m,1}$ if and only if it does not}
\\
\nonumber
&\text{ have subgraphs of the form in \eqref{eq:A}}
\end{align}

Let us denote by $F$ the set that consists of the graphs in 
\eqref{eq:A}. 
Since $B$ is equivalent to a tree $T$ of type $D_{m,1}$,
there exist a sequence of basic moves $\phi_{c_r,a_r},...,\phi_{c_2,a_2}, \phi_{c_1,a_1}$
such that  $B=\phi_{c_r,a_r}\circ...\circ\phi_{c_2,a_2}\circ\phi_{c_1,a_1}(T)$. We will show by 
induction on $r$ that 
each graph $B_i=\phi_{c_i,a_i}...\phi_{c_2,a_2}\phi_{c_1,a_1}(T)$, $1 \leq i \leq r$
, in particular $B=B_r$, contains a subgraph that belongs to $F$. 
The basis of the induction is the fact that
any tree of type $D_{m,1}$ does not contain any
subgraph that belongs to $F$. The induction follows from the following lemma:

\begin{lemma}
\label{lem:rec-A}
If a graph $G$ contains a subgraph $X \in F$, then, for any basic move $\phi_{c,a}$, the graph $\phi_{c,a}(G)$ contains a subgraph $X'$ which belongs to $F$.

\end{lemma}

Since the basic move exchanges $c$ by $c+a$ and fixes the other
elements of $G$, we may assume that $c$ is in X.
It follows from a direct check that
\begin{align} 
\label{eq:A-star} 
&\text{if $a \in X$, then $\phi_{c,a}(G)$ contains a subgraph $X'$ which belongs to $F$}.
\end{align}

If $a \notin X$ and it is 
not connected to any vertex $v$ in $X$ such that $v \ne c$, then the graph $\phi_{c,a}(X)=X$ is in $F$. 
Let us now assume that $a\notin X$ and it is connected to a vertex $v\ne c$ in $X$. 
By \eqref{eq:A-star}, we may also assume that
$a$ is not contained in any subgraph which is in $F$.
Then the subgraph $Xa$ induced by $X \cup \{a\}$ is 
of the form in Fig.~\ref{fig:Ahat-0} as could be verified easily.
We note that the graph $X'=\phi_{c,a}(Xa)$ belongs to $F$.

\begin{figure}[ht]

\setlength{\unitlength}{1.8pt}

\begin{center}

\begin{picture}(120,80)(-110,-30)

\put(-50,10){\line(1,-3){10}}
\put(-50,10){\line(-1,-3){10}}
\put(-50,10){\circle*{2.0}}
\put(-47,13){\makebox(0,0){$a$}}

\thicklines



\put(-20,20){\circle*{2.0}}
\put(-20,0){\circle*{2.0}} 
\put(-20,20){\circle*{2.0}}
\put(-40,40){\circle*{2.0}}
\put(-60,40){\circle*{2.0}}
\put(-80,20){\circle*{2.0}}
\put(-60,-20){\circle*{2.0}}
\put(-80,0){\circle*{2.0}}
\put(-40,-20){\circle*{2.0}}

\put(-20,0){\line(-1,-1){20}}
\put(-20,0){\line(0,1){20}}
\put(-20,20){\line(-1,1){20}}
\put(-40,40){\line(-1,0){20}}
\put(-60,40){\line(-1,-1){20}}
\put(-80,20){\line(0,-1){20}}
\put(-80,0){\line(1,-1){20}}
\put(-60,-20){\line(1,0){20}}

\put(-19,-4){\makebox(0,0){$x_s$}}
\put(-40,43){\makebox(0,0){$x_i$}}
\put(-84,0){\makebox(0,0){$x_1$}}
\put(-84,20){\makebox(0,0){$x_2$}}
\put(-40,-23){\makebox(0,0){$c$}}
\put(-60,-23){\makebox(0,0){$v$}}
\put(1,-4){\makebox(0,0){$s \geq 2$}}

\end{picture}
\end{center}

\caption{} 

\label{fig:Ahat-0}

\end{figure}

We complete the proof of Theorem~\ref{th:subE6} by the following lemma.

\begin{lemma}

\label{lem:subE6}

Let $G$ be a connected graph that contains a subgraph $X$ 
which is equivalent to $E_6$. Then, for any $a,c$ in $G$
such that $\Omega(a,c)=1$,
the graph $\phi_{c,a}(G)$ contains a subgraph $X'$ which is equivalent to the 
Dynkin tree $E_6$.

\end{lemma}

If $c \notin X$, then we may take $X'=X$. 
If $c \in X$ and $a \in X$, then we may take $X'=\phi_{c,a}(X)$. 
Let us now assume that $c \in X$ and $a \notin X$. 
We write $Xc$ to denote the (connected) graph induced by the
vertices in $X \cup \{c\}$ and denote by $U$ the linear span of $Xc$. We will show that 

\begin{align} 
\label{eq:E6} 
&\text{the graph $\phi_{c,a}(Xc)$ contains a subgraph $X'$ which is equivalent to $E_6$.} 
\end{align} 

By Theorem~\ref{th:Janssen classification}, the graph $Xc$
is equivalent to a tree of the form in  Fig.~\ref{fig:JanssenA} or Fig.~\ref{fig:JanssenC}.
Let us first assume that $Xc$, hence $Y=\phi_{c,a}(Xc)$, is equivalent to 
a tree of the form in Fig.~\ref{fig:JanssenA}. 
Then $dim (U_{00})=1$ by Proposition~\ref{th:complement of U000}.
Let $y$ in $U_{00}$ be such that $y \ne 0$, and let $b$ be a vertex in $Gr(Y,y)$
such that $X'=Y-\{y\}$ is connected. Then $Arf(Q_{X'})=1$, so $X'$ is equivalent
to $E_6$ by Theorem~\ref{th:Janssen classification}.

Let us now assume that $Xc$ is equivalent to a tree of the form in Fig.~\ref{fig:JanssenC}. 
We will show, using a case by case analysis, that the graph 
$Y=\phi_{c,a}(Xc)$ contains a subgraph $X'$ which 
is equivalent to $E_6$, and this will complete the proof of Lemma~\ref{lem:subE6}. 

For the remaing part of the proof, we assume that $Y=\{b_1,b_2,...,b_6,b_7\}$. We denote by $U$ the linear span of $Y$.  
We note that $dim(U_0)=1$ and $dim(U_{00})=dim(U_{000})=0$.

\credit{Case 1} 
\emph{$Y$ is a cycle of length 7}

Any cycle $C$ of length $r\geq 4$ is equivalent to a tree of type $D_{r-2,2}$. Since $Xc$, hence $Y$, is equivalent to a tree that contains $E_6$, this case is not possible by Theorem~\ref{pr:broom}.

\credit{Case 2} 
\emph{$Y$ contains a cycle $C$ of length 6.}

Let us assume without loss of generality that $C=[b_1,b_2,...,b_6]$.
Since $dim(U_{00})=0$, the vertex $b_7$ is connected to an odd number 
of vertices in $C$ (otherwise the vector $b_1+b_2+...+b_6$ is in $U_{00}$). 
Thus $Y$ is one of the graphs in Figures~\ref{fig:E7-21}-\ref{fig:E7-23}
each contains a subgraph equivalent to $E_6$ which is marked by thick lines.




\begin{figure}[ht]

\setlength{\unitlength}{1.8pt}

\begin{center}

\begin{picture}(70,50)(10,-25)

\put(54,0){\makebox(0,0){$b_7$}}
\put(50,0){\circle*{2.0}}
\put(60,-20){\line(-1,0){30}}
\put(80,0){\line(-1,-1){20}}

\thicklines

\put(60,20){\circle*{2.0}}
\put(30,20){\circle*{2.0}}
\put(30,-20){\circle*{2.0}}
\put(60,-20){\circle*{2.0}}
\put(80,0){\circle*{2.0}}
\put(10,0){\circle*{2.0}}

\put(27,-21){\makebox(0,0)}
\put(26,23){\makebox(0,0)}

\put(10,0){\line(1,1){20}}
\put(30,-20){\line(-1,1){20}}
\put(30,20){\line(1,0){30}}
\put(60,20){\line(1,-1){20}}
\put(30,20){\line(1,-1){20}}

\end{picture}

\end{center}

\caption{} 

\label{fig:E7-21}

\end{figure}






\begin{figure}[ht]

\setlength{\unitlength}{1.8pt}

\begin{center}

\begin{picture}(170,50)(-40,-25)

\put(54,0){\makebox(0,0){$b_7$}}
\put(50,0){\circle*{2.0}}
\put(80,0){\line(-1,-1){20}}
\put(60,20){\line(1,-1){20}}

\thicklines

\put(60,20){\circle*{2.0}}
\put(30,20){\circle*{2.0}}
\put(30,-20){\circle*{2.0}}
\put(60,-20){\circle*{2.0}}
\put(80,0){\circle*{2.0}}
\put(10,0){\circle*{2.0}}

\put(27,-21){\makebox(0,0)}
\put(26,23){\makebox(0,0)}

\put(10,0){\line(1,1){20}}
\put(30,-20){\line(-1,1){20}}
\put(30,20){\line(1,0){30}}
\put(30,20){\line(1,-1){20}}
\put(30,-20){\line(1,1){20}}
\put(10,0){\line(1,0){40}}
\put(60,-20){\line(-1,0){30}}

\end{picture}

\end{center}

\caption{} 

\label{fig:E7-221}

\end{figure}



\begin{figure}[ht]

\setlength{\unitlength}{1.8pt}

\begin{center}

\begin{picture}(170,50)(70,-25)

\put(164,0){\makebox(0,0){$b_7$}}
\put(160,0){\circle*{2.0}}
\put(190,0){\line(-1,-1){20}}
\put(170,20){\line(1,-1){20}}

\thicklines

\put(170,20){\circle*{2.0}}
\put(140,20){\circle*{2.0}}
\put(140,-20){\circle*{2.0}}
\put(170,-20){\circle*{2.0}}
\put(190,0){\circle*{2.0}}
\put(120,0){\circle*{2.0}}

\put(27,-21){\makebox(0,0)}
\put(26,23){\makebox(0,0)}

\put(120,0){\line(1,1){20}}
\put(140,-20){\line(-1,1){20}}
\put(140,20){\line(1,0){30}}
\put(140,20){\line(1,-1){20}}
\put(140,-20){\line(1,1){20}}
\put(170,-20){\line(-1,2){10}}
\put(170,-20){\line(-1,0){30}}

\end{picture}

\end{center}

\caption{} 

\label{fig:E7-222}

\end{figure}



\begin{figure}[ht]

\setlength{\unitlength}{1.8pt}

\begin{center}

\begin{picture}(170,50)(-40,-25)

\put(54,0){\makebox(0,0){$b_7$}}
\put(50,0){\circle*{2.0}}
\put(80,0){\line(-1,-1){20}}
\put(60,20){\line(1,-1){20}}

\thicklines

\put(60,20){\circle*{2.0}}
\put(30,20){\circle*{2.0}}
\put(30,-20){\circle*{2.0}}
\put(60,-20){\circle*{2.0}}
\put(80,0){\circle*{2.0}}
\put(10,0){\circle*{2.0}}

\put(27,-21){\makebox(0,0)}
\put(26,23){\makebox(0,0)}

\put(10,0){\line(1,1){20}}
\put(30,-20){\line(-1,1){20}}
\put(30,20){\line(1,0){30}}
\put(60,20){\line(-1,-2){10}}
\put(60,-20){\line(-1,2){10}}
\put(10,0){\line(1,0){40}}
\put(60,-20){\line(-1,0){30}}

\end{picture}

\end{center}

\caption{} 

\label{fig:E7-223}

\end{figure}



\begin{figure}[ht]

\setlength{\unitlength}{1.8pt}

\begin{center}

\begin{picture}(170,50)(-40,-25)

\put(54,0){\makebox(0,0){$b_7$}}
\put(50,0){\circle*{2.0}}
\put(80,0){\line(-1,-1){20}}
\put(60,20){\line(1,-1){20}}

\thicklines

\put(60,20){\circle*{2.0}}
\put(30,20){\circle*{2.0}}
\put(30,-20){\circle*{2.0}}
\put(60,-20){\circle*{2.0}}
\put(80,0){\circle*{2.0}}
\put(10,0){\circle*{2.0}}

\put(27,-21){\makebox(0,0)}
\put(26,23){\makebox(0,0)}

\put(10,0){\line(1,1){20}}
\put(30,-20){\line(-1,1){20}}
\put(30,20){\line(1,0){30}}
\put(30,20){\line(1,-1){20}}
\put(30,-20){\line(1,1){20}}
\put(10,0){\line(1,0){40}}
\put(60,-20){\line(-1,0){30}}
\put(60,20){\line(-1,-2){10}}
\put(60,-20){\line(-1,2){10}}

\end{picture}

\end{center}

\caption{} 

\label{fig:E7-23}

\end{figure}

\credit{Case 3} 
\emph{$Y$ contains a cycle $C$ of length 5.}

Let us assume without loss of generality that $C=[b_1,b_2,...,b_5]$.
Since $dim(U_{00})=0$, there is a vertex, say $b_6$, in $Y$ that is connected to an 
odd number of vertices in $C$ (otherwise the vector $b_1+b_2+...+b_5$ is in $U_{00}$). 
Then the graph induced by $\{b_1,b_2,...,b_6\}$ is one of the
graphs in Fig.~\ref{fig:E7-3a} or Fig.~\ref{fig:E7-3b}
each equivalent to $E_6$. 

\begin{figure}[ht]

\setlength{\unitlength}{1.8pt}

\begin{center}

\begin{picture}(140,50)(0,-25)



\thicklines


\put(-10,0){\circle*{2.0}}
\put(30,20){\circle*{2.0}}
\put(20,-20){\circle*{2.0}}
\put(40,-20){\circle*{2.0}}
\put(50,0){\circle*{2.0}}
\put(10,0){\circle*{2.0}}

\put(130,20){\circle*{2.0}}
\put(120,-20){\circle*{2.0}}
\put(140,-20){\circle*{2.0}}
\put(150,0){\circle*{2.0}}
\put(110,0){\circle*{2.0}}
\put(130,0){\circle*{2.0}}

\put(-11,-3){\makebox(0,0){$b_6$}}
\put(19,-23){\makebox(0,0)}
\put(8,4){\makebox(0,0){$b_2$}}

\put(10,0){\line(1,1){20}}
\put(10,0){\line(1,-2){10}}
\put(30,20){\line(1,-1){20}}
\put(50,0){\line(-1,-2){10}}
\put(20,-20){\line(1,0){20}}
\put(10,0){\line(-1,0){20}}

\put(110,0){\line(1,1){20}}
\put(110,0){\line(1,-2){10}}
\put(130,20){\line(1,-1){20}}
\put(150,0){\line(-1,-2){10}}
\put(120,-20){\line(1,0){20}}
\put(130,0){\line(-1,-2){10}}
\put(130,0){\line(-1,0){20}}
\put(130,0){\line(0,1){20}}

\end{picture}

\end{center}

\caption{} 

\label{fig:E7-3a}

\end{figure}

\begin{figure}[ht]

\setlength{\unitlength}{1.8pt}

\begin{center}

\begin{picture}(140,50)(0,-25)



\thicklines


\put(30,20){\circle*{2.0}}
\put(20,-20){\circle*{2.0}}
\put(40,-20){\circle*{2.0}}
\put(50,0){\circle*{2.0}}
\put(10,0){\circle*{2.0}}
\put(30,0){\circle*{2.0}}

\put(130,0){\circle*{2.0}}
\put(130,20){\circle*{2.0}}
\put(120,-20){\circle*{2.0}}
\put(140,-20){\circle*{2.0}}
\put(150,0){\circle*{2.0}}
\put(110,0){\circle*{2.0}}


\put(10,0){\line(1,1){20}}
\put(10,0){\line(1,-2){10}}
\put(30,20){\line(1,-1){20}}
\put(50,0){\line(-1,-2){10}}
\put(20,-20){\line(1,0){20}}
\put(30,0){\line(-1,-2){10}}
\put(30,0){\line(-1,0){20}}
\put(30,0){\line(1,0){20}}

\put(110,0){\line(1,1){20}}
\put(110,0){\line(1,-2){10}}
\put(130,20){\line(1,-1){20}}
\put(150,0){\line(-1,-2){10}}
\put(120,-20){\line(1,0){20}}
\put(130,0){\line(-1,-2){10}}
\put(130,0){\line(-1,0){20}}
\put(130,0){\line(0,1){20}}
\put(130,0){\line(1,0){20}}
\put(130,0){\line(1,-2){10}}

\end{picture}

\end{center}

\caption{} 

\label{fig:E7-3b}

\end{figure}

\credit{Case 4} 
\emph{$Y$ contains a cycle $C$ of length 4.}

Let us assume, without loss of generality, that $C=[b_1,b_2,b_3,b_4]$.
Since $dim(U_{00})=0$, there is a vertex $v$ in $\{b_5,b_6,b_7\}$ such that
$v$ is connected to an odd number of vertices in $C$.

\credit{Subcase 4.1}\emph{There is a vertex $v\notin C$  such that
$v$ is connected to precisely 3 vertices in $C$}.

Let us assume, without loss of generality, that $v$ 
is connected to the vertices $b_1,b_2,b_3$.
Since $dim(U_{00})=0$, there is a vertex $v'\ne v$ such that $v' \notin C$  
and connected to precisely one vertex, say $b_1$, in $\{b_1,b_3\}$ (otherwise $b_1+b_3$ in $U_{00}$).

\credit{Subsubcase 4.1.1} 
\emph{$v'$ is not connected to any vertex in $\{b_2,b_4\}$} 
Then the graph $X'$ induced by $C\cup \{v,v'\}$ is one of the graphs in 
Fig.~\ref{fig:E7-41}
and it is equivalent to $E_6$.

\begin{figure}[ht]

\setlength{\unitlength}{1.8pt}

\begin{center}

\begin{picture}(70,50)(-10,-10)



\thicklines

\put(-30,0){\circle*{2.0}}
\put(-10,0){\circle*{2.0}}
\put(-50,20){\circle*{2.0}}
\put(-30,20){\circle*{2.0}}
\put(-10,20){\circle*{2.0}}
\put(-20,10){\circle*{2.0}}

\put(70,0){\circle*{2.0}}
\put(90,0){\circle*{2.0}}
\put(70,20){\circle*{2.0}}
\put(90,20){\circle*{2.0}}
\put(110,40){\circle*{2.0}}
\put(70,40){\circle*{2.0}}

\put(10,-3){\makebox(0,0){}}
\put(110,-3){\makebox(0,0){}}
\put(-22,8){\makebox(0,0){$v$}}
\put(-50,23){\makebox(0,0){$v'$}}
\put(113,40){\makebox(0,0){$v$}}
\put(67,41){\makebox(0,0){$v'$}}
\put(-5,21){\makebox(0,0){$b_2$}}
\put(-31,-2){\makebox(0,0){$b_4$}}
\put(-31,24){\makebox(0,0){$b_1$}}
\put(-5,-2){\makebox(0,0){$b_3$}}

\put(-50,20){\line(1,0){20}}
\put(70,0){\line(1,0){20}}
\put(-30,0){\line(0,1){20}}
\put(-10,0){\line(0,1){20}}
\put(-30,20){\line(1,0){20}}
\put(-30,0){\line(1,0){20}}
\put(-10,0){\line(-1,1){20}}
\put(-10,20){\line(-1,-1){10}}

\put(70,0){\line(0,1){20}}
\put(90,0){\line(0,1){20}}
\put(70,20){\line(1,0){20}}
\put(70,40){\line(0,-1){20}}
\put(110,40){\line(-2,-1){40}}
\put(110,40){\line(-1,0){40}}
\put(110,40){\line(-1,-2){20}}
\put(110,40){\line(-1,-1){20}}

\end{picture}

\end{center}

\caption{} 

\label{fig:E7-41}

\end{figure}

\credit{Subsubcase 4.1.2} 
\emph{$v'$ is connected to precisely one vertex in $\{b_2,b_4\}$}

Let $X'$ be the graph induced by $C\cup \{v,v'\}$. Then we have the following:

(i) If $v'$ is not connected to $v$ then the graph $X'$ is one of the graphs in Fig.\ref{fig:E7-412} each equivalent to $E_6$.

(ii) If $v'$ is connected to $v$ then (i) can be applied
to $\phi_{v',b_1}(X')$ and conclude that $X'$ is equivalent to $E_6$.

\begin{figure}[ht]

\setlength{\unitlength}{1.8pt}

\begin{center}

\begin{picture}(70,50)(-10,-10)



\thicklines

\put(-30,0){\circle*{2.0}}
\put(-10,0){\circle*{2.0}}
\put(-50,20){\circle*{2.0}}
\put(-30,20){\circle*{2.0}}
\put(-10,20){\circle*{2.0}}
\put(-20,10){\circle*{2.0}}

\put(70,0){\circle*{2.0}}
\put(90,0){\circle*{2.0}}
\put(70,20){\circle*{2.0}}
\put(90,20){\circle*{2.0}}
\put(80,10){\circle*{2.0}}
\put(70,40){\circle*{2.0}}

\put(10,-3){\makebox(0,0){}}
\put(110,-3){\makebox(0,0){}}
\put(-22,8){\makebox(0,0){$v$}}
\put(-50,23){\makebox(0,0){$v'$}}
\put(78,8){\makebox(0,0){$v$}}
\put(67,41){\makebox(0,0){$v'$}}
\put(-5,21){\makebox(0,0){$b_2$}}
\put(-31,-2){\makebox(0,0){$b_4$}}
\put(-31,24){\makebox(0,0){$b_1$}}
\put(-5,-2){\makebox(0,0){$b_3$}}

\put(-50,20){\line(1,0){20}}
\put(-50,20){\line(1,-1){20}}
\put(70,0){\line(1,0){20}}
\put(-30,0){\line(0,1){20}}
\put(-10,0){\line(0,1){20}}
\put(-30,20){\line(1,0){20}}
\put(-30,0){\line(1,0){20}}
\put(-10,0){\line(-1,1){20}}
\put(-10,20){\line(-1,-1){10}}

\put(70,40){\line(1,-1){20}}
\put(70,0){\line(0,1){20}}
\put(90,0){\line(0,1){20}}
\put(70,20){\line(1,0){20}}
\put(70,40){\line(0,-1){20}}
\put(80,10){\line(1,1){10}}
\put(80,10){\line(-1,1){10}}
\put(80,10){\line(1,-1){10}}

\end{picture}

\end{center}

\caption{} 

\label{fig:E7-412}

\end{figure}

\credit{Subsubcase 4.1.3} 
\emph{$v'$ is connected to precisely two vertices in $\{b_2,b_4\}$}

Let $X'$ be the graph induced by $C\cup \{v,v'\}$.
Then Subcase 4.1.1 applies to $\phi_{v',b_1}(X')$ so $X'$ is equivalent to $E_6$.

\credit{Subcase 4.2} 
\emph{There is no vertex $v \in \{b_5,b_6,b_7\}$ such that
$v$ is connected to 3 vertices in $C$.}

Since $dim(U_{00})=0$, there are vertices 
$v,v' \notin C$ which are connected
to adjacent vertices in $C$.
Then the graph $X'$ induced by $C\cup\{v,v'\}$ is 
one of the graphs in Fig.~\ref{fig:E7-42a} each equivalent to $E_6$.

\begin{figure}[ht]

\setlength{\unitlength}{1.8pt}

\begin{center}

\begin{picture}(70,50)(-10,-10)



\thicklines

\put(-50,0){\circle*{2.0}}
\put(-10,0){\circle*{2.0}}
\put(-30,20){\circle*{2.0}}
\put(-30,0){\circle*{2.0}}
\put(-50,20){\circle*{2.0}}
\put(-30,20){\circle*{2.0}}
\put(-10,20){\circle*{2.0}}

\put(30,0){\circle*{2.0}}
\put(70,0){\circle*{2.0}}
\put(50,20){\circle*{2.0}}
\put(50,0){\circle*{2.0}}
\put(30,20){\circle*{2.0}}
\put(50,20){\circle*{2.0}}
\put(70,20){\circle*{2.0}}

\put(10,-3){\makebox(0,0){}}
\put(110,-3){\makebox(0,0){}}

\put(-50,20){\line(1,0){20}}
\put(-30,0){\line(0,1){20}}
\put(-10,0){\line(0,1){20}}
\put(-30,20){\line(1,0){20}}
\put(-50,0){\line(1,0){40}}

\put(30,20){\line(1,0){20}}
\put(50,0){\line(0,1){20}}
\put(70,0){\line(0,1){20}}
\put(50,20){\line(1,0){20}}
\put(30,0){\line(1,0){40}}
\put(30,0){\line(0,1){20}}

\end{picture}

\end{center}

\caption{} 

\label{fig:E7-42a}

\end{figure}

\credit{Case 5} 
\emph{$Y$ contains two adjacent triangles sharing precisely one
common edge.}

Let $T=\{b_1,b_2,b_3,b_4\}$ be the subgraph formed by the adjacent triangles
sharing the common edge $[b_2,b_3]$. We note that the graph $T'=\phi_{b_1,b_2}(T)$
is a cycle of length 4. By our analysis in Case 4, the graph $T'$ is contained in a
graph $E$ which is equivalent to $E_6$. Then $T$ is contained in 
$X'=\phi_{b_1,b_2}(E)$.

\credit{Case 6} 
\emph{Y contains a subgraph $D$ of the form $D_{2,2}$.}

Let $D=\{a_1,a_2,c_1,c_2\}$ be indexed as in Fig.~\ref{fig:broom}.
We note that the graph $D'=\phi_{a_1,a_2}(D)$
is induced by two adjacent triangles. 
By our analysis in Case 5, the graph $D'$ is contained in a
graph $E$ which is equivalent to $E_6$. Then $D$ is contained in 
$X'=\phi_{a_1,a_2}(E)$. 

\credit{Case 7} 
\emph{None of the above cases happen.}

Then Y is equivalent to a tree of type $D_{m,1}$ by \eqref{eq:A2}, which contradicts to our assumption that $Y$ is equivalent to a tree that contains $E_6$ Proposition~\ref{pr:E_6 ne D}. This completes the proof of Theorem~\ref{th:subE6}.

\subsection{Proof of Theorem~\ref{th:D-V000}}
\label{subsec:D-V000}
Let us first assume that $X$ has precisely 4 vertices. 
Then it is equivalent to a tree $X'$ of type $D_{2,2}$.  For $X'$, Theorem~\ref{th:D-V000} follows from Lemma~\ref{lem:A-d}.
Let us now assume that $X$ has at least 5 vertices, i.e. $X=[c_1,....,c_r]$ is a cycle 
whose length $r$ is grater than or equal to $5$. We claim that $x=c_1+c_2+...+c_r$ is in $V_{000}$. 
Let us note that $Gr(B',x)$ is a cycle of length $r-1$, where $B'=\phi_{b_{r-1},b_r}(B)$.
Then the claim follows from an induction.

The remaining part of Theorem~\ref{th:D-V000} follows from Lemma~\ref{lem:A-d}.

\subsection{Proof of Theorem~\ref{th:D-V000-char}}
\label{subsec:D-V000-char}
We will show that, for any $B$ which is equivalent to a tree that contains $E_6$,
there is a subgraph $X$ of the form in $\eqref{eq:A}$ such that $\FF_2^X\cap V_{000}= \{0\}$.
For such a $B$, there is a subgraph $E$ which is equivalent to $E_6$. Let us denote
by $U$ the linear span of vectors contained in $E$. Then $dim(U_0)=0$ by Proposition~\ref{th:complement of U000}. On the other hand,
since $E$ is not equivalent to any tree of type $D_{m,1}$, it contains a subgraph $X$ such that $X$ is
of the form in \eqref{eq:A} (Theorem~\ref{th:subE6}). Then $\FF_2^X\cap V_{000}= \{0\}$ and we are done.






\subsection{Proof of Theorem~\ref{th:D-V000-span}}
\label{subsec:D-V000-span}
For any basis $B$ which is equivalent to a tree of type $D_{m,k}$ we denote
by $B_{000}$ the set that consists of vectors $x\in V_{000}$ such that $Gr(B,x)$ is 
contained in a subgraph of the form in \eqref{eq:A}. Then Theorem~\ref{th:D-V000-span}
is equivalent to the following statement:
\begin{align} 
\label{eq:D-th} 
&\text{For any basis $B$ which is equivalent to a tree of type $D_{m,k}$, 
the set $B_{000}$ spans $V_{000}.$}
\end{align}
 
Let us first note that the theorem holds for $B$ which
is a tree of type $D_{m,k}$ (Proposition~\ref{th:complement of U000}). 
Thus, to prove \eqref{eq:D-th}, it is enough to prove the following statement:  
 
\begin{align} 
\label{eq:D} 
&\text{if $x$ in $B_{000}$ then, for any basis
$B'=\phi_{c,a}(B)$, the vector $x$ is in the linear span of $B'_{000}$.}
\end{align}

It follows from a direct check that
\begin{align} 
\label{eq:D-star} 
&\text{if $a,c$ are in $Gr(B,x)$, then \eqref{eq:D} is satisfied.}
\end{align}

Let us now assume, without loss of generality, that 
$c \in Gr(B,x)$, and $a \notin Gr(B,x)$. We will establish
\eqref{eq:D} using a case by case analysis.

\credit{Case 1} \emph{$Gr(B,x)$ has precisely two vertices.}

In this case $Gr(B,x)$ is contained in a subgraph $X$ which is a tree of type $D_{2,2}$ or a cycle of length 4. 
Let $x=c+b$ be the expansion of $x$ in the basis $B$. We note
that $a$ is connected to both $c$ and $b$ because $x\in V_0$. 
Then $Gr(B',x)$ is the triangle induced by $\{c+a,a,b\}$ and
it is in $B'_{000}$. 
 


\credit{Case 2} \emph{$Gr(B,x)$ has precisely three vertices} 

We note that in this case $Gr(B,x)$ is a triangle.
Let $x=c+b_1+b_2$ be the expansion of $x$ in the basis $B$.
Then $Gr(B',x)$ is the cycle induced by $\{c+a,b_1,b_2,a\}$.
In particular it is in $B'_{000}$.



\credit{Case 3}\emph{$Gr(B,x)$ has precisely 4 vertices}

We note that in this case $Gr(B,x)$ is a cycle of length 4.
Let us assume that $x=c+b_1+b_2+b_3$ is the expansion of $x$ in the basis $B$.

\credit{Subcase 3.1}\emph{$a$ is connected to precisely one
vertex, say b, in $\{b_1,b_2,b_3\}$.}
If $b$ is connected to $c$, then $Gr(B',x)$ is a cycle of length of 5 
so $x\in B'_{000}$.
Let us now assume that $b$ is not connected to $c$. 
We may also assume that $b=b_2$ (see Fig.~\ref{fig:D-V000-3}(i)). 
Then the vectors $x_1=b_1+a$ and $x_2=b_2+b_3+c+a$ are in $V_{000}$
by Theorem~\ref{th:D-V000}. They are also in the linear span of $B'_{000}$ by \eqref{eq:D-star},  
so is $x=x_1+x_2$.

\credit{Subcase 3.2.}\emph{$a$ is connected to all of the vertices in $\{b_1,b_2,b_3\}$.}
In this case, the vectors $x_1=c+a+b_1$ and $x_2=b_2+a+b_3$ are in $V_{000}$ 
by Theorem~\ref{th:D-V000} (see Fig.~\ref{fig:D-V000-3}(ii)).
They are also in the linear span of $B'_{000}$ by \eqref{eq:D-star},  
so is $x=x_1+x_2$.

\begin{figure}[ht]
\setlength{\unitlength}{1.8pt}
\begin{center}
\begin{picture}(70,45)(-35,-5)
\thicklines
\put(-60,0){\circle*{2.0}}
\put(-30,0){\circle*{2.0}}
\put(-60,30){\circle*{2.0}}
\put(-30,30){\circle*{2.0}}
\put(-45,15){\circle*{2.0}}
\put(35,15){\circle*{2.0}}
\put(-60,0){\line(1,0){30}}
\put(-60,30){\line(1,0){30}}
\put(-60,0){\line(0,1){30}}
\put(-30,0){\line(0,1){30}}
\put(-60,0){\line(1,1){30}}
\put(-63,-2){\makebox(0,0){$b_2$}}
\put(-63,33){\makebox(0,0){$b_3$}}
\put(-28,32){\makebox(0,0){$c$}}
\put(-27,-3){\makebox(0,0){$b_1$}}
\put(-43,13){\makebox(0,0){$a$}}
\put(-43,-15){\makebox(0,0){(i)}}

\put(20,0){\circle*{2.0}}
\put(50,0){\circle*{2.0}}
\put(20,30){\circle*{2.0}}
\put(20,0){\circle*{2.0}}
\put(50,30){\circle*{2.0}}
\put(20,30){\line(1,0){30}}
\put(20,0){\line(1,0){30}}
\put(20,0){\line(0,1){30}}
\put(50,0){\line(0,1){30}}
\put(20,0){\line(1,1){30}}
\put(50,0){\line(-1,1){30}}

\put(18,-2){\makebox(0,0){$b_2$}}
\put(18,34){\makebox(0,0){$b_3$}}
\put(53,32){\makebox(0,0){$c$}}
\put(52,-3){\makebox(0,0){$b_1$}}
\put(38,15){\makebox(0,0){$a$}}
\put(38,-15){\makebox(0,0){(ii)}}


\end{picture}

\end{center}

\caption{}

\label{fig:D-V000-3}

\end{figure}

\credit{Case 4} \emph{$Gr(B,x)$ contains precisely 5 vertices.}

In this case $Gr(B,x)=[c,b_1,b_2,b_3,b_4]$ is a cycle of length 5. 

\credit{Subcase 4.1} \emph{$a$ is connected to precisely one
vertex, say $b$, in $\{b_1,b_2,b_3,b_4\}$.}

If $b$ is connected to $c$, then $Gr(B',x)$ is a cycle of length 5,
so $x\in B'_{000}$. 
Let us now assume that $b$ is not connected to $c$.
We may also assume, without loss of generality, that $b=b_2$ (see Fig.~\ref{fig:D-V000-4}(i)).
Then the vectors $x_1=b_1+a, x_2=c+a+b_2+b_3+b_4$ are in $V_{000}$ by Theorem.
They are also in the linear span of $B'_{000}$ by \eqref{eq:D-star}, so is $x=x_1+x_2$.

\credit{Subcase 4.2} \emph{$a$ is connected to precisely three
vertices in $\{b_1,b_2,b_3,b_4\}$.}

In this case, we may assume that $Gr(B,x)$ is as in 
Fig.~\ref{fig:D-V000-4}(ii). 
We note that $x_1=c+b_1+a$ and $x_2=a=b_2+b_3+b_4$ are in $V_{000}$. 
They are also in the linear span of $B'_{000}$ by \eqref{eq:D-star}, so is $x=x_1+x_2$.

\begin{figure}[ht]

\setlength{\unitlength}{1.8pt}

\begin{center}

\begin{picture}(140,50)(0,-25)



\thicklines


\put(30,0){\circle*{2.0}}
\put(30,20){\circle*{2.0}}
\put(20,-20){\circle*{2.0}}
\put(40,-20){\circle*{2.0}}
\put(50,0){\circle*{2.0}}
\put(10,0){\circle*{2.0}}

\put(130,20){\circle*{2.0}}
\put(120,-20){\circle*{2.0}}
\put(140,-20){\circle*{2.0}}
\put(150,0){\circle*{2.0}}
\put(110,0){\circle*{2.0}}
\put(130,0){\circle*{2.0}}

\put(31,-33){\makebox(0,0){(i)}}
\put(31,-3){\makebox(0,0){$a$}}
\put(32,23){\makebox(0,0){$b_1$}}
\put(19,-23){\makebox(0,0){$b_4$}}
\put(8,3){\makebox(0,0){$c$}}
\put(53,3){\makebox(0,0){$b_2$}}
\put(41,-23){\makebox(0,0){$b_3$}}

\put(132,-33){\makebox(0,0){(ii)}}
\put(131,23){\makebox(0,0){$b_1$}}
\put(132,-3){\makebox(0,0){$a$}}
\put(119,-23){\makebox(0,0){$b_4$}}
\put(108,3){\makebox(0,0){$c$}}
\put(153,3){\makebox(0,0){$b_2$}}
\put(141,-23){\makebox(0,0){$b_3$}}

\put(10,0){\line(1,1){20}}
\put(10,0){\line(1,-2){10}}
\put(30,20){\line(1,-1){20}}
\put(50,0){\line(-1,-2){10}}
\put(20,-20){\line(1,0){20}}
\put(10,0){\line(1,0){40}}

\put(110,0){\line(1,1){20}}
\put(110,0){\line(1,-2){10}}
\put(130,20){\line(1,-1){20}}
\put(150,0){\line(-1,-2){10}}
\put(120,-20){\line(1,0){20}}
\put(130,0){\line(-1,-2){10}}
\put(110,0){\line(1,0){40}}
\put(130,0){\line(0,1){20}}

\end{picture}

\end{center}

\caption{} 

\label{fig:D-V000-4}

\end{figure}

\credit{Case 5}{$Gr(B,x)$ has at least 6 vertices.}

We note that the graph $Gr(B,x)=[c,b_1,...,b_r]$ is a cycle of length $r, r\geq 6$. 

\credit{Subcase 5.1}{$a$ is connected to precisely one
vertex, say $b$ in ${b_1,b_2,b_3,b_4,b_5}$.}

If $b$ is connected to $c$, then
$Gr(B',x)$ is a cycle of length of $r+1$, so $x\in B'_{000}$ 
Let us now assume that $b$ is not connected to $c$
and denote the graph induced by $Gr(B,x) \cup \{a\}$ by $Xa$.
One could easily check the following:
if $Xa$ does not contain a cycle of length 4 then it 
contains a subgraph which is $E_6$, 
contradicting to our assumption that $B$ does not contain a subgraph equivalent to $E_6$ . 
If $Xa$ contains a cycle of lenth 4, then
it is of the form in Fig~\ref{fig:D-V000-51}.
By the same arguments as in Case 5, one could show that $x$ in the
linear span of $B'_{000}$.

\credit{Subcase 5.2}
Let us now assume that $a$ is connected to (an odd number) $k$
vertices in ${b_1,b_2,...,b_r}$ with $k\geq 1$.
One could easily check that $Xa$ 
contains a subgraph equivalent to $E_6$ unless it is 
of the form in Fig~\ref{fig:D-V000-52}. 
By the same arguments as in Case 5, one could show that $x$ in the
linear span of $B'_{000}$.


\begin{figure}[ht]
\setlength{\unitlength}{1.8pt}
\begin{center}
\begin{picture}(170,50)(70,-25)

\thicklines

\put(163,0){\makebox(0,0){$a$}}
\put(160,0){\circle*{2.0}}
\put(170,-20){\line(2,1){10}}
\put(170,20){\line(2,-1){10}}

\put(170,20){\circle*{2.0}}
\put(140,20){\circle*{2.0}}
\put(140,-20){\circle*{2.0}}
\put(170,-20){\circle*{2.0}}
\put(190,0){\circle*{2.0}}
\put(185,10){\circle*{2.0}}
\put(185,-10){\circle*{2.0}}

\put(120,0){\circle*{2.0}}

\put(27,-21){\makebox(0,0)}
\put(26,23){\makebox(0,0)}

\put(120,0){\line(1,1){20}}
\put(140,-20){\line(-1,1){20}}
\put(140,20){\line(1,0){30}}
\put(140,20){\line(1,-1){20}}
\put(140,-20){\line(1,1){20}}
\put(170,-20){\line(-1,0){30}}

\end{picture}
\end{center}
\caption{} 
\label{fig:D-V000-51}
\end{figure}

\begin{figure}[ht]

\setlength{\unitlength}{1.8pt}

\begin{center}

\begin{picture}(170,50)(70,-25)

\thicklines

\put(163,0){\makebox(0,0){$a$}}
\put(160,0){\circle*{2.0}}
\put(170,-20){\line(2,1){10}}
\put(170,20){\line(2,-1){10}}

\put(170,20){\circle*{2.0}}
\put(140,20){\circle*{2.0}}
\put(140,-20){\circle*{2.0}}
\put(170,-20){\circle*{2.0}}
\put(190,0){\circle*{2.0}}
\put(185,10){\circle*{2.0}}
\put(185,-10){\circle*{2.0}}

\put(120,0){\circle*{2.0}}
\put(125,10){\circle*{2.0}}
\put(125,-10){\circle*{2.0}}

\put(27,-21){\makebox(0,0)}
\put(26,23){\makebox(0,0)}

\put(140,20){\line(-2,-1){10}}
\put(140,-20){\line(-2,1){10}}
\put(140,20){\line(1,0){30}}
\put(140,20){\line(1,-1){20}}
\put(140,-20){\line(1,1){20}}
\put(170,-20){\line(-1,2){10}}
\put(170,20){\line(-1,-2){10}}
\put(170,-20){\line(-1,0){30}}

\end{picture}
\end{center}
\caption{} 
\label{fig:D-V000-52}
\end{figure}

\subsection{Proof of Theorem~\ref{th:xxx}}
\label{subsec:orbit-D}



Let us first introduce some notation. For any vector $v\in V$, we 
denote by $C(B,v)$ the set of connected components of $Gr(B,v)$.
We denote by $M(B,v)$ the set of maximal complete subgraphs (of $Gr(B,v)$)
with at least three vertices.
By some abuse of notation, we will also use $\bar{x}$
to denote the vector that corresponds to the graph $\bar{x}$.
We will denote by $d_B(\bar{x})$ the expression on the right side of \eqref{eq:xxx},
thus we will show that $d(x)=d_B(\bar{x})$ for any $x\in V-V_0$ and any $B$ which is equivalent 
to a tree of type $D_{m,k}$.

Since $\bar{x}$ is \emph{maximal}, each connected component of $Gr(B,\bar{x})$ is
equivalent to a tree of type $D_{m,1}$ (Theorems~\ref{th:D-V000-span}, \ref{th:subE6}).
Let $A=\{a_1,...,a_r\}, r\geq 3$ be a maximal complete subgraph which 
is contained in a connected component $C$ of $\bar{x}$.
For any vertex $b$ in $C$ such that $b$ is not in $A$, 
we denote by $A(b)$ the vertex in $A$ which is closest to $b$ 
(such a vertex exists by \eqref{eq:A}). 
For any vertex $a\in A$, we define $H(a)=\{b \in C-A: A(b)=a\} \cup \{a\}$
and for a pair $\{a_i,a_j\}$ of vertices  in $A$, we define 
$H(a_i+a_j)=\{b \in C-A: A(b)=a_i \mathrm{\:or\:} A(b)=a_j\} \cup \{a_i+a_j\}$.
If $r$ is even (resp. odd)
let $$B'=\phi_{a_{r-1},a_r}\circ...\circ \phi_{a_3,a_4}\circ \phi_{a_1,a_2} 
(\mathrm{resp. \:\:} B'=\phi_{a_{r-2},a_{r-1}}\circ...\circ \phi_{a_3,a_4}\circ \phi_{a_1,a_2}).$$
Then 
$C(B',\bar{x})=(C(B,\bar{x})-C) \cup \{H(a_1+a_2),H(a_3+a_4),...,H(a_{r-1}+a_r)\}$ 
(resp. 
$C(B',\bar{x})=(C(B,\bar{x})-C) \cup \{H(a_1+a_2),H(a_3+a_4),...,H(a_{r-2}+a_{r-1})\}, H(a_r)\}$),
and $M(B',x)=M(B,x)-A$. We note that $d_B(\bar{x})=d_{B'}(\bar{x})$ and
and $\bar{x}$ is minimal. Continuing this procedure, we obtain
a basis $B''$ which is equivalent $B$ such that each connected component of $Gr(B'',\bar{x})$
is a chain. Then $d_B(\bar{x})=d_{B''}(\bar{x})$ is equal to the number of connected components of $Gr(B'',\bar{x})$.
(Note that $\bar{x}$ is also minimal with respect to $B'$).
Since each connected component of $Gr(B'',\bar{x})$ corresponds to a vector in $\Delta$ (c.f. Subsection~\ref{sec:delta decomposition}), we have $d_{B''}(\bar{x})=d(\bar{x})=d(x)$ by Lemma~\ref{lem:A-d}, Theorems~\ref{pr:broom}, 
\ref{th:delta decomposition} and we are done.

\section{Orbits of groups generated by symplectic
transvections of a linearly independent subset}

\label{sec:independent}

In this section, we will prove Theorem~\ref{non-trivial orbits} after some preliminary statements. Throughout the section, $B$ denotes a linearly independent
subset which is not a basis in a finite dimensional $\FF_2$-space $V$ 
equipped with
the alternating form $\Omega$. We always assume that $Gr(B)$ is connected.
We denote by $U$ the linear span of $B$. We note that each $\gamma \in 
\Gamma_B$ preserves cosets in $V/U$, so we only need to describe   
$\Gamma_B$-orbits in each coset $v+U$. If $v+U=U$,
then our previous results apply, so we will always
consider the action of $\Gamma_B$ on a coset $v+U\ne U$.
We note that the set $B\cup \{v\}$ is linearly independent
and there is the associated graph $Gr(B\cup \{v\})$ as defined in Section~\ref{sec:background}.
We denote by $V^{\Gamma_B}$ the set of 
vectors in $V$ which are fixed by $\Gamma_B$. 
As before, $U_0$ denotes the kernel of the form $\Omega \vert_U$ and
$\Delta$ is the $\Gamma_B$-orbit that contains $B$ 
(\cite[Proposition~3.1]{BHI}). The spaces $U_{00}$ and $U_{000}$ are defined as in Section~\ref{sec:background}.

If $(v+U)\cap V^{\Gamma_B}$ is non-empty, then our previuous results
allows one to describe all $\Gamma_B$ orbits in $v+U$. More precisely,
we have the following statement.

\begin{proposition}
\label{pr:orbits type D2}

Suppose that $(v+U)\cap V^{\Gamma_B}$ is non-empty and contains a vector $v+u$. Then $\Gamma_B$-orbits in $v+U$ are parallel translates of $\Gamma_B$-orbits in $U$ by $v+u$.
\end{proposition}



\begin{figure}[ht]
\setlength{\unitlength}{1.8pt}
\begin{center}
\begin{picture}(70,20)(-10,-10)
\thicklines

\put(-10,0){\circle*{2.0}}
\put(0,0){\circle*{2.0}}
\put(10,0){\circle*{2.0}}
\put(20,0){\circle*{2.0}}
\put(30,0){\circle*{2.0}}
\put(40,0){\circle*{2.0}}
\put(50,0){\circle*{2.0}}
\put(60,0){\circle*{2.0}}
\put(60,10){\circle*{2.0}}
\put(60,-10){\circle*{2.0}}
\put(20,10){\circle*{2.0}$^v$}

\put(-10,0){\line(1,0){70}}
\put(50,0){\line(1,1){10}}
\put(50,0){\line(1,-1){10}}
\put(10,0){\line(1,1){10}}
\put(20,0){\line(0,1){10}}
\put(40,0){\line(-2,1){20}}
\put(20,10){\line(1,0){40}}

\end{picture}
\end{center}
\caption{}
\label{fig:6}
\end{figure}

Our next result gives a sufficient condition for $(v+U)\cap V^{\Gamma_B}$ to be empty.

\begin{proposition}
\label{pr:U000 and trivial}

If $\Omega(v,U_{000}) \ne \{0\}$, then $(v+U)\cap V^{\Gamma_B}=\emptyset$.

\end{proposition}

\proof

We may assume, by Theorems~\ref{pr:broom},~\ref{th:Janssen classification},  
that $B$ is equivalent to one of the four trees indexed as in Fig.~\ref{fig:broom}, 
Fig.~\ref{fig:JanssenA},
Fig.~\ref{fig:JanssenB} and Fig.~\ref{fig:JanssenC}.

Let us introduce the numbers $r$ and $t$ as follows: 
if $B$ is equivalent to the tree in Fig.~\ref{fig:broom} with $m=2$, 
then $r=m$, $t=k+1$ and we set $c_t=a_1$; 
if $B$ is equivalent to th tree in Fig.~\ref{fig:broom} with $m>2$, then $r=m$, $t=k$; 
if $B$ is equivalent to the tree in Fig.~\ref{fig:JanssenA} or Fig.~\ref{fig:JanssenB}; 
then $r=2n-1$, $t=p+1$; 
if $B$ is equivalent to the tree in Fig.~\ref{fig:JanssenC}, 
then $r=2n$, $t=p$. 
 
If $\Omega(v,U_{000}) \ne \{0\}$, then the set $I=\{c_i: \Omega(v,c_i)=1\}$ is a 
non-empty, proper subset of $\{c_1,...,c_t\}$ by Proposition~\ref{th:complement of U000}.
We assume, without loss of generality, that $I=\{c_1,...,c_s\}, s<t$. 
Let $x \in U$. If $x$ contains $a_r$, then $\tau_{c_{s+1}}(v+x) \ne v+x$. If $x$  
does not contain $a_m$, then $\tau_{c_1}(v+x) \ne v+x$. Thus
$v+x \notin V^{\Gamma_B}$ for any $x\in U$, i.e. 
$(v+U)\cap V^{\Gamma_B} =\emptyset$. 
\endproof


\subsection{Proof of Theorem~\ref{non-trivial orbits}(i)}


We first note that there exist at least two $\Gamma_B$-orbits in $v+U$
because, for any $u\in \Delta$ such that $\Omega(u,v)=0$, we have
$Q_{B\cup\{v\}}(v+u)=0$, so $v$ and $v+u$ lie in different orbits;
here the existence of $u$ follows from our assumption that $dim(U) \geq 2$.




Since $B$ does not contain any subgraph which is
equivalent to $E_6$, it is equivalent to a tree of type
$D_{m,k}$ with $m\geq 2, k\geq 1$ (Theorem~\ref{th:subE6}). Since basic moves preserve
the associated qudratic forms (c.f. Section~\ref{sec:background}), we may take $B$ as in Fig.~\ref{fig:broom} 
with the same indexing.
A typical graph of $B\cup {v}$ is given in Fig.~\ref{fig:6}.
As a first step, we will disconnect $v$ from $a_i$'s using basic moves.

\begin{lemma}
\label{}
There exists $\alpha \in \Gamma_B$ such that $\Omega(\alpha(v),a_j)=0$ for $j=1,...,m$.
\end{lemma}

If $m=2$, we set $c_{k+1}=a_1$ for convenience. 
Since $\Omega(v,U_{000})\ne \{0\}$, 
we may assume that $c_1\in I=\{c_i: \Omega(v,c_i)=1\}$
by Proposition~\ref{th:complement of U000}.

For any $w\in v+U$, we define $A(w)=\{a_i:\Omega(a_i,w)=1\}$.
If $A(w) \ne \emptyset$, we let $i(w)=max\{i:a_i\in A(w)\}$. 
Let us write $\alpha_i=\tau_{a_{i+1}}...\tau_{a_m}\tau_{c_1}$ for $i<m$.
If $i=i(w)<m$, then 
$A(\alpha_i(w))=(A(w)-\{a_i\})\cup \{a_{i+1}\}$, so we have
$A(\alpha_{m-1}...\alpha_{i+1}\alpha_i(w))=(A(w)-\{a_i\})\cup \{a_{m}\}$. We also note that 
if $i(w)=m$ and $\Omega(w,c_1)=1$ then $A(\tau_{c_1}(w))=A(w)-\{a_m\}$.
Thus, by induction on $i(v)$ if necessary, we obtain
$\alpha \in \Gamma_B$ such that $A(\gamma(v))$ is a proper 
subset of $A(v)$.
By induction on the cardinality of $A(v)$, we obtain $\gamma\in \Gamma_B$
such that $A(\gamma (v))=\emptyset$. This completes the proof of the lemma.

For the remaining part of the proof of Theorem~\ref{non-trivial orbits}(i)
we assume, without loss of generality, that 
$I=\{b\in B: \Omega(b,v)=1\}=\{c_1,...,c_s\}$, where $s<k$ if $m>2$ and $s \leq k+1$
if $m=2$ (here $c_{k+1}=a_1$). 


\begin{lemma}
\label{}
Let $f$ be the linear map on the span of $B\cup \{v\}$ defined as follows:
$f(c_i)=c_1$ for $c_i$ in $I$, 
and $f(b)=b$ for $b\in B\cup\{v\}$ such that $b\notin I$. Then,
for any $z \in U$, the vectors $v+z$ and $f(v+z)$ lie in the same $\Gamma_B$-orbit.
\end{lemma}

If $c_i$ is not contained in $v+z$ for any $i\in I$, then 
$f(v+z)=v+z$ and we are done. Let us assume that $v+z$ contains 
$c_1,...,c_l$ from $I$. 
If $v+z$ does not contain $a_m$ and $l$ is odd (resp. even), 
then $f(v+z)=\tau_{c_l}...\tau_{c_2}(v+z)$ (resp. $f(v+z)=\tau_{c_l}...\tau_{c_1}(v+z)$.
Let us now assume that $v+z$ contains $a_m$. 
If $\Omega(a_m,z)=1$ and $l$ is odd (resp. even), 
then $f(v+z)=\tau_{a_m}\tau_{c_2}...\tau_{c_l}\tau_{a_m}(v+z)$
(resp. $f(v+z)=\tau_{a_m}\tau_{c_1}...\tau_{c_l}\tau_{a_m}(v+z)$).
If $\Omega(a_m,z)=0$ and $l$ is odd (resp. even), then
$f(v+z)=\tau_{c_{s+1}}\tau_{a_m}\tau_{c_2}...\tau_{c_l}\tau_{a_m}\tau_{c_{s+1}}(v+z)$
(resp. $f(v+z)=\tau_{c_{s+1}}\tau_{a_m}\tau_{c_1}...\tau_{c_l}\tau_{a_m}\tau_{c_{s+1}}(v+z)$). This completes the proof of the lemma

Thus Theorem~\ref{non-trivial orbits}(i) is equivalent to the following statement:

$f(v+x)$ and $f(v+y)$ lie in the same $\Gamma_B$-orbit if
and only if $Q_{B\cup\{v\}}(f(v+x))=Q_{B\cup\{v\}}(f(v+y))$.

We complete the proof by the following lemma. We recall
that $Q_B(\alpha(w))=Q_B(w)$ for any $w\in V$ and $\alpha \in \Gamma_B$.

\begin{lemma}

\label{}

For any $u\in U$, the vector $f(v+u)$ is in the orbit of either $v$ or 
$v+a_m$.

\end{lemma}

We first note that $v$ and $v+a_m$ lie in different $\Gamma_B$-orbits
because $Q_{B\cup\{v\}}(v)=1$ and $Q_{B\cup\{v\}}(v+a_m)=0$.  
To prove the lemma, we will first show that there is $\gamma \in \Gamma_B$
such that $\gamma(f(v+u))$ is contained in the chain $A$
formed by $a_1,...,a_m,c_1,v$. Recall that $f(v+u)$
does not contain any of $\{c_2,...,c_s\}$. 
Let us first assume that $f(v+u)$ has a component $u_m$ that contains $a_m$. 
If $u_m$ contains vertices $c_{j_1},...,c_{j_l} \subset \{c_{s+1},...,c_k\}$, 
then for $\gamma=\tau_{c_{j_1}}...\tau_{c_{j_l}}$, $\gamma (f(v+u))$ is contained in 
$A$. Now let us assume that $u$ does not contain $a_m$. If $\Omega(a_m,f(v+u))=1$, then
$\tau_{a_m}(f(v+u))=f(v+u)+a_m$ contains $a_m$ and the previous arguments
apply. If $\Omega(a_m,f(v+u))=0$, then 
$\tau_{a_m}\tau_{c_1}(f(v+u))=f(v+u)+c_1+a_m$ contains $a_m$
and we apply the previous arguments.

Thus, we may assume that $f(v+u)=v+f(u)$ is contained in the chain formed by $a_1,...,a_m,c_1,v$.
We may also assume that $\Omega(v,f(u))=0$, (otherwise we can write 
$f(v+u)=c_1+a_m+a_{m-1}+...+a_{m-i}+x$ where $x\in span(a_1,...,a_i)$ and consider
$\tau_{c_1}\tau_{a_m}\tau_{a_{m-1}}...\tau_{a_{m-i}}(f(v+u))$.
Then, by Theorem~\ref{pr:broom}, there exists 
$\beta \in \Gamma_{\{a_1,...,a_m\}}$ such that 
$\beta(f(v+u))=v+a_m+a_{m-2}+...+a_{m-2r}$ for some $r \geq 0$.
If $r=0$, then we are done. We suppose $r\geq 1$. Then, for
$\alpha=\tau_{c_1}\tau_{a_m}\tau_{a_{m-1}}\tau_{c_{s+1}}\tau_{a_m}\tau_{c_1}\tau_{a_{m-2}}\tau_{a_{m-1}}\tau_{a_m}
\tau_{c_{s+1}}$, we have 
$\alpha(v+a_m+a_{m-2}+...+a_{m-2r})=
v+a_{m-4}+...+a_{m-2r}$, which has two less components than 
$\beta (f(v+u))=v+a_m+a_{m-2}+...+a_{m-2r}$. Continuing
this process, we will have $f(v+x)$ in the orbit of either $v$ or $v+a_m$, which proves the 
lemma. This also completes the proof of Theorem~\ref{non-trivial orbits}(i). 



\subsection{Proof of Theorem~\ref{non-trivial orbits}(ii)}

If $(v+U)\cap V^{\Gamma_B}\ne \emptyset$, then the statement follows from Proposition~\ref{pr:U000 and trivial}. 
Let us now assume that $(v+U)\cap V^{\Gamma_B}\ne \emptyset$.
By Theorem~\ref{th:subE6}, we may take $B$ as in Fig.~\ref{fig:broom} with the same indexing.
If $\Omega(v,U_{000})= \{0\}$, then  $v$ is connected to none of the $c_i$'s or
connected to all of them as in Fig.~\ref{fig:12a} and Fig.~\ref{fig:12b}. If $v$ is connected to all of $c_i$'s, then
$v+a_m$ is connected to only $a_i$'s, so we may assume that $v$ is connected only to $a_i$'s.
Suppose $v$ is connected to $a_i$, i.e. $\Omega(v,a_i)=1$, but $\Omega(v,a_j)=0$ for
$j=i+1,..,m$. Then $v+a_{i-1}$ will not be connected to $a_i,...,a_m$. Continuing this way,
we will have a $w=v+u$ connected to only $a_1$. Thus, $B\cup\{w\}$ will be of type $D_{m+1,k}$.

Let us now prove that,
\begin{align} 
\label{eq:ii} 
&\text{for $x,y \in U$, the vectors $w+x$ and $w+y$ lie in the same 
$\Gamma_B$-orbit}
\\
\nonumber 
&\text{if and only if they lie in the same $\Gamma_{B\cup\{w\}}$-orbit.}
\end{align}

The "only if" part follows from the fact that
$\Gamma_B$ is a subgroup of $\Gamma_{B\cup\{w\}}$. To prove the "if" part, 
we assume, by Theorem~\ref{pr:broom}, that $f(w+x)$ and $f(w+y)$ have
the same number of connected components. Then it is easy to see the following:
there is $\alpha , \beta \in \Gamma_B$ such that $\alpha(w+x)=w+x'$, 
$\beta(w+y)=w+y'$ with $\Omega(w,x')$=0 and $\Omega(w,y')$=0. We will
show that $w+x'$ and $w+y'$ lie in the same $\Gamma_B$-orbit. We note that
$x',y'\in S=span(a_2,...,a_m,c_1,...,c_k)$ and $f(x')$ and $f(y')$ have the 
same number of connected components). Thus there exists $\gamma \in \Gamma_S$
such that $\gamma(x')=y'$. Since $\Omega(w,s)=0$ for any $s\in S$, 
we have $\gamma(w+x')=w+y'$ and we are done.    
\endproof

\begin{figure}[ht]
\setlength{\unitlength}{1.8pt}
\begin{center}
\begin{picture}(70,20)(-10,-10)

\thicklines
\put(0,0){\circle*{2.0}}
\put(10,0){\circle*{2.0}}
\put(70,0){\circle*{2.0}$^v$}
\put(20,0){\circle*{2.0}}
\put(30,0){\circle*{2.0}}
\put(40,0){\circle*{2.0}}
\put(50,0){\circle*{2.0}}
\put(60,0){\circle*{2.0}}
\put(60,10){\circle*{2.0}}
\put(60,-10){\circle*{2.0}}

\put(-10,0){\line(1,0){80}}
\put(50,0){\line(1,1){10}}
\put(50,0){\line(1,-1){10}}
\put(70,0){\line(-1,1){10}}
\put(70,0){\line(-1,-1){10}}

\end{picture}
\end{center}
\caption{}
\label{fig:12a}
\end{figure}

\begin{figure}[ht]
\setlength{\unitlength}{1.8pt}
\begin{center}
\begin{picture}(70,20)(-10,-10)
\thicklines

\put(0,0){\circle*{2.0}}
\put(10,0){\circle*{2.0}}
\put(20,0){\circle*{2.0}}
\put(30,0){\circle*{2.0}}
\put(40,0){\circle*{2.0}}
\put(50,0){\circle*{2.0}}
\put(60,0){\circle*{2.0}}
\put(60,10){\circle*{2.0}}
\put(60,-10){\circle*{2.0}}
\put(20,10){\circle*{2.0}$^v$}

\put(-10,0){\line(1,0){70}}
\put(50,0){\line(1,1){10}}
\put(50,0){\line(1,-1){10}}
\put(10,0){\line(1,1){10}}
\put(20,0){\line(0,1){10}}
\put(40,0){\line(-2,1){20}}

\end{picture}
\end{center}
\caption{}
\label{fig:12b}
\end{figure}

\begin{example}

\label{ex}

\rm{

For this example, let $B=\{b_1,b_2,b_3,b_4,b_5,b_6\}$ and let $V$ denote the
vector space over $\FF_2$ with basis $B\cup \{v_1,v_2,v_3,v_4\}$.
We denote by $\Omega$ the skew-symmetric form given in Fig.~\ref{fig:ex}.
As above, $U$ denotes the linear span of $B$ in $V$. 
We determine $\Gamma_B$-orbits in $V$ as follows:

It is easy to see that $c=b_2+b_4+b_6 \in U_{000}$. 
Let $V_{1,4}$ denote the linear span of the set $\{v_1,v_4\}$.
Since $\Omega(v_2,c)=\Omega(v_3,c)=1$ and $\Omega(v_1,c)=\Omega(v_4,c)=0$,
by Theorem~\ref{non-trivial orbits}, we have the following:

If $v\in v_2+V_{1,4}$ or $v\in v_3+V_{1,4}$, 
then the $\Gamma_B$-orbits in the coset $v+U$
are intersections of $v+U$ with the sets $Q_{B\cup\{v\}}^{-1}(0)$ and $Q_{B\cup\{v\}}^{-1}(1)$
(so there are 16 of such orbits).

The remaining $\Gamma$-orbits are contained in cosets $v+U$
where $v\in S=span(\{v_2+v_3,v_1,v_4\})$.

To proceed, we first notice that $v_2+v_3+b_5$ and $v_1+v_4+b_6+b_5$ 
are fixed by $\Gamma_B$. Now we note the following fact:
for $v,w \in V$, if $v+w$ is fixed by $\Gamma_B$, then $\Gamma_B$-orbits
in $v+U$ are parallel translates of $\Gamma_B$-orbits in $w+U$ by $v+w$
(because $\alpha(v+u)=\alpha(v+w+w+u)=v+w+\alpha(w+u)$ 
for all $\alpha\in \Gamma_B, u\in U$). We also note that 
the one-element $\Gamma_B$-orbits in $U$ are the vectors 
$\{0,b_2+b_6+b_4, b_1+b_5+b_3+b_4, b_1+b_5+b_3+b_2+b_6\}$. 
The non-trivial $\Gamma_B$-orbits have representatives
$b_5$ and $b_5+b_6$. Also all $\Gamma_B$-orbits in $v_1+U$ are non-trivial and they have
representatives $v_1,v_1+b_5,v_1+b_5+b_6$.

Thus, the total number of $\Gamma_B$-orbits is $16+4\cdot 6+4\cdot 3=16+24+12=52$.

According to \cite{SSVZ}, there is a bijection between $\Gamma_B$-orbits and
connected components of the reduced double Bruhat cell $L^{w_0,e}(\RR)$ 
for $W=S_5$ and $$w_0=s_1s_3s_2s_4s_1s_3s_2s_4s_1s_3$$
where $s_i=(i,i+1)$ are adjacent transpositions. We note that $w_0$ is the longest 
element of $W$ and $\ii=(1,3,2,4,1,3,2,4,1,3)$ is its reduced word and
the set of bounded indices $B(\ii)$(see \cite{SSVZ}) is $B$.
We remark that the total number of $\Gamma_B$-orbits 
agrees with the result given in \cite{SSV}.

}

\end{example}

\begin{figure}[ht]

\setlength{\unitlength}{1.8pt}

\begin{center}

\begin{picture}(100,50)(0,0)

\thicklines

\put(0,0){\circle*{2.0}$^{v_1}$}

\put(10,20){\circle*{2.0}$^{v_2}$}

\put(20,10){\circle*{2.0}$^{v_3}$}

\put(30,30){\circle*{2.0}$^{v_4}$}

\put(40,0){\circle*{2.0}$^{b_1}$}

\put(50,20){\circle*{2.0}$^{b_2}$}

\put(60,10){\circle*{2.0}$^{b_3}$}

\put(70,30){\circle*{2.0}$^{b_4}$}

\put(80,0){\circle*{2.0}$^{b_5}$}

\put(90,20){\circle*{2.0}$^{b_6}$}

\put(0,0){\line(1,0){80}}

\put(20,10){\line(1,0){40}}

\put(10,20){\line(1,0){80}}

\put(30,30){\line(1,0){40}}

\put(40,0){\line(2,1){20}}

\put(40,0){\line(-2,1){20}}

\put(50,20){\line(2,1){20}}

\put(50,20){\line(1,-1){10}}

\put(50,20){\line(-2,1){20}}

\put(50,20){\line(-3,-1){30}}

\put(60,10){\line(2,-1){20}}

\put(60,10){\line(3,1){30}}

\put(70,30){\line(2,-1){20}}

\end{picture}

\end{center}

\caption{}

\label{fig:ex}

\end{figure}

\section{Orbits of Groups generated by non-symplectic tranvections}

\label{sec:5}

In this section, we prove Theorem~\ref{th:GSV} and give an example. 
One could easily note that Theorem~\ref{th:GSV} follows from the following lemma
which extends \cite[Theorem~3]{GSV} to an arbitrary bilinear form.

\begin{lemma}
\label{th:arbitrary}

Let $V$ be an $\FF_2$ space equipped with a non-skew-symmetric bilinear form $\Omega(u,v)$. Suppose that $B$ is a linearly independent subset of $V$ such that $\Omega(b,b)=0$ for all $b\in B$ and $Gr(B)$ is connected.
Let $x$ be a vector in $V$ and let $B_L$ be an arbitrary connected subgraph of such that the group $\Gamma_{B_L}$ does not fix $x$. Suppose that $b \in B-B_L$ and let $P=[b_0 \in T(x),...,b_k=b]$ be a shortest path that connects $b$ to $T(x)={b \in B_L:\tau_b(x) \ne x}$. If $\Omega \vert_{P_L(b,x)}$ is not alternating and $\Omega(b_i,b_{i+1})=1$ for $i=0,...,k-1$,
then the the vector $x+b$ lies in the $\Gamma_B$ orbit that contains $x$. 

\end{lemma}

We prove the lemma by modifying the proof of \cite[Theorem~3]{GSV}.
Let us first note that $\Omega(b_i,b_j)=0$ for $0\leq i<i+1<j\leq k$ because
$P$ is a shortest path.
Let $l=\min\{i\geq 1:\Omega(b_{i-1},b_i)=1 \mathrm{\:but\:} \Omega(b_i,b_{i-1})=0\}$.


We will prove the lemma by induction on $k\geq 0$, the length of $P$. The case $k=0$ is
clear: $\tau_b(x)=x+b$. 
Suppose the statement of the lemma holds for $P$ of length less than $k$.
Let us first assume that $\tau_{b_j}(x)=x$ for
$j=l,...,k$. Then $\tau_{b_k}...\tau_{b_0}(x)=x+b_0+...+b_{l-1}+b_l+...+b_k$.
Since $\tau_{b_{l-1}}(x+b_0+...+b_k)=(x+b_0+...+b_k)+b_{l-1}$, we have
$$\tau_{b_k}...\tau_{b_l}\tau_{b_0}...\tau_{b_{l-1}}\tau_{b_k}...\tau_{b_0}(x)=x+b.$$
Let us now assume that there exists $l-1<j<k+1$ such that $\tau_{b_j}(x)\ne x$ and let $m$ be the
one closest to $b_k=b$. We note that $z=\tau_{b_k}...\tau_{b_m}(x)=x+b_m+...+b_k$ is
in the same orbit as $x$.
Then,
the length of the shortest path connecting $T(z)$ to $b_{k-1}$ is less than $k$, hence by
the induction hypothesis, $z+b_{k-1}=x+b_m+...+b_{k-2}$ lies in the orbit of $x$.
Applying the same procedure to $b_{k-2},...,b_m$, we will have $x+b_k$ in the orbit of $x$.
\endproof

\begin{example}

\label{exx}

\rm{

Let $V$ be the vector space over $\FF_2$
with basis $B'=\{b_1,...,b_6,v_1,v_2,v_3\}$. Let us introduce the following
sets: $B=\{b_1,b_2,b_3,b_4,b_5,b_6\}$, $B_R=\{b_1,b_2,b_4,b_5\}$, 
$B_L=\{b_3,b_6\}$, $C_R=\{v_1,v_2\}$, $C_L=\{v_3\}$, $R=B_R\cup C_R$
$L=B_L\cup C_L$. 
We denote by $\Omega$ the bilinear form on $V$ defined by $Gr(B')$ given in Fig.~\ref{fig:exx}
such that $\Omega \vert_{\FF_2^R}$ and $\Omega \vert_{\FF_2^L}$ 
are alternating and $\Omega(a_r,a_l)=0$ for any $a_r \in R$ and $a_l \in L$.   

To determine the $\Gamma_B$-orbits in $V$, we first note that 
the vectors $\{b_1+b_2+b_4, b_1+b_5, b_2+b_4+b_5\}\subset (\FF _2^{B_U})$.
Let us write $V_{1,2}=$linear span of $\{v_1,v_2,v_3+b_6\}$, which is
the set of fixed points of $\Gamma_{B_L}$.
Then for any $v\in V_{1,2}$ such that $v\notin \{0, v_3+v_1+b_6\}$, 
we have $\Omega(v,(\FF _2^{B_R})_{000})\ne 0$.
Thus, by Proposition~\ref{th:GSV} and
Theorem~\ref{non-trivial orbits}, we have the following:
For any $v\in V_{1,2}$ such that $v\notin \{0,v_3+v_1+b_6\}$, the $\Gamma_B$-orbits in the coset $v+\FF _2^{B_R}$
are intersections of $v+\FF _2^{B_R}$ with the sets $Q_{B\cup\{v\}}^{-1}(0)$ and $Q_{B\cup\{v\}}^{-1}(1)$
(so there are 12 such orbits).
We also note that $v_3+v_1+b_6+b_4$ is fixed by $\Gamma_{B_R}$, 
thus the $\Gamma_{B_R}$(hence $\Gamma_B$)-orbits in $v_3+v_1+b_6+\FF _2^{B_R}$ are parallel
translates of $\Gamma_{B_R}$-orbits in $\FF _2^{B_R}$ by Proposition~\ref{pr:orbits type D2}.
The one-element $\Gamma_{B_R}$-orbits in $\FF _2^{B_R}$ are $0, b_2+b_4+b_5, b_2+b_4+b_1, b_1+b_5$
and the remaining one is represented by $b_1$, so there is a total of 10 $\Gamma_B$-orbits 
in $v+\FF _2^{B_R}$
for any $v\in \{0,v_3+v_1+b_6\}$. According to Theorem~\ref{th:GSV}, the remaining orbits are 
represented by the vectors
$v_3,v_3+v_1,v_3+v_2,v_3+v_1+v_2,b_3,b_3+v_1,b_3+v_2,b_3+v_1+v_2$.

Thus, the total number of $\Gamma_B$-orbits is 12+10+8=30.

According to \cite{Z}, $\Gamma_B$-orbits are in bijection with
the connected components of the reduced double Bruhat cell $L^{w_0,e}(\RR)$ 
for $W$ of type $B_3$, where $w_0$ is te longest element in its Weyl group. 
We remark that the total number of $\Gamma_B$-orbits 
agrees with the result obtained in \cite{GSV}.

}

\end{example}

\begin{figure}[ht]

\setlength{\unitlength}{1.8pt}

\begin{center}

\begin{picture}(70,40)(-10,0)

\thicklines

\put(0,0){\circle*{2.0}$^{v_1}$}

\put(10,10){\circle*{2.0}$^{v_2}$}

\put(20,20){\circle*{2.0}$^{v_3}$}

\put(30,0){\circle*{2.0}$\: \: ^{b_1}$}

\put(40,10){\circle*{2.0}$\: ^{b_2}$}

\put(50,20){\circle*{2.0}$^{b_3}$}

\put(60,0){\circle*{2.0}$\: \: ^{b_4}$}

\put(70,10){\circle*{2.0}$\: ^{b_5}$}

\put(80,20){\circle*{2.0}$^{b_6}$}

\put(0,0){\line(1,0){60}}

\put(10,10){\line(1,0){60}}

\put(20,20){\line(1,0){60}}

\put(30,0){\line(-2,1){20}}

\put(30,0){\line(1,1){10}}

\put(40,10){\line(2,-1){20}}

\put(40,10){\line(-2,1){20}}

\put(40,10){\line(1,1){10}}

\put(70,10){\line(-2,1){20}}

\put(60,0){\line(1,1){10}}

\put(70,10){\line(1,1){10}}

\end{picture}

\end{center}

\caption{}

\label{fig:exx}

\end{figure}

\section{Acknowledgements}

I am grateful to my graduate advisor Andrei Zelevinsky for his support and many suggestions.
I also thank Misha Kogan and Josh Scott for helpful discussions.

\end{document}